\documentclass[a4paper,12pt]{article}

\pdfoutput=1

\usepackage[paper=letterpaper,margin=1in]{geometry}

\usepackage{amsmath,amssymb,amsfonts,epsfig,cite,setspace,bigstrut,longtable,array,breqn,url,color}

\usepackage{hyperref}

\allowdisplaybreaks

\usepackage{lipsum}

\let\OLDthebibliography\thebibliography
\renewcommand\thebibliography[1]{
  \OLDthebibliography{#1}
  \setlength{\parskip}{0pt}
  \setlength{\itemsep}{0pt plus 0ex}
}

\definecolor{purple}{cmyk}{0,0.8,0,0.4}
\definecolor{db}{cmyk}{1,0.1,0.2,0.6}
\definecolor{dg}{cmyk}{1,0,1,0.7}
\definecolor{bl}{cmyk}{1,0.7,0.5,0.2}
\definecolor{yl}{rgb}{0.2,0.7,0.2}
\definecolor{bl2}{cmyk}{0.7,0.4,0,0.5}
\definecolor{red2}{cmyk}{0,1,1,0.8}
\definecolor{red3}{cmyk}{0,0.7,1,0.7}
\definecolor{gr2}{cmyk}{1,0.2,0.7,0.6}
\definecolor{nb}{rgb}{0.1,0.1,0.5}
\definecolor{ng}{rgb}{0,0.8,0}
\definecolor{brown}{rgb}{0.6,0.3,0.2}
\definecolor{newred}{cmyk}{0,1,1,1}

\newcommand{\red}{\color{red}}
\newcommand{\blue}{\color{blue}}

\newcommand{\purple}{\color{purple}}

\begin{document}
\newcommand{\todo}[1]{{\bf ?????!!!! #1 ?????!!!!}\marginpar{$\Longleftarrow$}}
\newcommand{\nn}{\nonumber}
\newcommand{\tr}{\mathop{\rm Tr}}
\newcommand{\ch}{\rm Ch}
\newcommand{\comment}[1]{}

\newcommand{\cM}{{\cal M}}
\newcommand{\cW}{{\cal W}}
\newcommand{\cN}{{\cal N}}
\newcommand{\cH}{{\cal H}}
\newcommand{\cK}{{\cal K}}
\newcommand{\cZ}{{\cal Z}}
\newcommand{\cO}{{\cal O}}
\newcommand{\cB}{{\cal B}}
\newcommand{\cC}{{\cal C}}
\newcommand{\cD}{{\cal D}}
\newcommand{\cE}{{\cal E}}
\newcommand{\cF}{{\cal F}}
\newcommand{\cR}{{\cal R}}
\newcommand{\IA}{\mathbb{A}}
\newcommand{\IB}{\mathbb{B}}
\newcommand{\IP}{\mathbb{P}}
\newcommand{\IQ}{\mathbb{Q}}
\newcommand{\IH}{\mathbb{H}}
\newcommand{\IR}{\mathbb{R}}
\newcommand{\IC}{\mathbb{C}}
\newcommand{\IF}{\mathbb{F}}
\newcommand{\IM}{\mathbb{M}}
\newcommand{\II}{\mathbb{I}}
\newcommand{\IZ}{\mathbb{Z}}
\newcommand{\re}{{\rm Re}}
\newcommand{\im}{{\rm Im}}
\newcommand{\sym}{{\rm Sym}}

\newcommand{\tmat}[1]{{\tiny \left(\begin{matrix} #1 \end{matrix}\right)}}
\newcommand{\mat}[1]{\left(\begin{matrix} #1 \end{matrix}\right)}
\newcommand{\diff}[2]{\frac{\partial #1}{\partial #2}}
\newcommand{\gen}[1]{\langle #1 \rangle}
\newcommand{\ket}[1]{| #1 \rangle}
\newcommand{\jacobi}[2]{\left(\frac{#1}{#2}\right)}

\def\acts{\curvearrowright}

\newcommand{\drawsquare}[2]{\hbox{%
\rule{#2pt}{#1pt}\hskip-#2pt
\rule{#1pt}{#2pt}\hskip-#1pt
\rule[#1pt]{#1pt}{#2pt}}\rule[#1pt]{#2pt}{#2pt}\hskip-#2pt
\rule{#2pt}{#1pt}}
\newcommand{\fund}{\raisebox{-.5pt}{\drawsquare{6.5}{0.4}}}
\newcommand{\antifund}{\overline{\fund}}

\newtheorem{theorem}{\bf THEOREM}
\def\thetheorem{\thesection.\arabic{theorem}}
\newtheorem{proposition}{\bf PROPOSITION}
\def\thetheorem{\thesection.\arabic{proposition}}
\newtheorem{observation}{\bf OBSERVATION}
\def\thetheorem{\thesection.\arabic{observation}}

\def\theequation{\thesection.\arabic{equation}}
\newcommand{\setall}{\setcounter{equation}{0}
        \setcounter{theorem}{0}}
\newcommand{\setequation}{\setcounter{equation}{0}}

~\\
\vskip 1cm

\begin{center}
{\Large \bf Sporadic and Exceptional}
\end{center}
\medskip

\vspace{.4cm}
\centerline{
{\large Yang-Hui He}$^1$ \&
{\large John McKay}$^2$
}
\vspace*{3.0ex}

\begin{center}
{\it
{\small
{${}^{1}$ 
Department of Mathematics, City University, London, EC1V 0HB, UK and \\
Merton College, University of Oxford, OX14JD, UK and\\
School of Physics, NanKai University, Tianjin, 300071, P.R.~China \\
\qquad
{\rm \url{hey@maths.ox.ac.uk}}\\
}
\vspace*{1.5ex}
{${}^{2}$ 
Department of Mathematics and Statistics,\\
Concordia University, 1455 de Maisonneuve Blvd.~West,\\
Montreal, Quebec, H3G 1M8, Canada\\
\qquad 
{\rm \url{mckay@encs.concordia.ca}}
}
}}
\end{center}

\begin{abstract}
We study the web of correspondences linking the exceptional Lie algebras $E_{8,7,6}$ and the sporadic simple groups Monster, Baby and the largest Fischer group.
This is done via the investigation of classical enumerative problems on del Pezzo surfaces in relation to the cusps of certain subgroups of $PSL(2,\IR)$ for the relevant McKay-Thompson series in Generalized Moonshine. 
We also study Conway's sporadic group, as well as its association with the Horrocks-Mumford bundle.
\end{abstract}

\newpage

\tableofcontents

~\\
~\\
~\\
~\\

\section{Introduction and Summary}\setall
The classification of mathematical structures oftentimes confronts the dichotomy between the {\it regular} and the {\it exceptional} wherein the former generically organizes into some infinite families whilst the latter tantalizes with what at first may seem an eclectic collage but out of whose initial disparity emerges striking order.
The earliest and perhaps most well-known example is that of the classification of symmetries in $\IR^3$, the regulars are the infinite families of the cyclic and dihedral groups of the regular $n$-gon and the exceptionals are the symmetry groups of the five Platonic solids.
Similarly, the classification of Lie algebras gave us the Dynkin diagrams of the infinite families of classical groups which are the familiar isometries of vector spaces, in addition, there are the five exceptional diagrams.
Indeed, the relation between the $\IR^3$ symmetries and the simply laced cases of the (affine) Lie algebras has come to be known as the McKay Correspondence, which now occupies a cornerstone of modern algebraic geometry and representation theory.

Another highlight example is of course the classification of finite simple groups which, after decades of arduous work, is by now complete.
The regulars here are the infinite families of cyclic groups of prime order, as well as the classical Lie groups over finite fields, while the exceptionals are known as the 26 sporadic groups, the largest of which is the Monster, of tremendous size.
Here, the second author's old observation that $196,884 = 196,883 + 1$, relating the Monster and the elliptic $j$-function, prompted the field of Moonshine \cite{CN,FLM}.

These above illustrative cases indeed exemplify how astute recognition of exceptional structures can help to unravel new mathematics of profound depth.
 Of particular curiosity is the less well-known fact - in parallel to the above identity - that the constant term of the $j$-invariant, viz., $744$, satisfies $744 = 3 \times 248$. The number 248 is, of course, the dimension of the adjoint of the largest exceptional algebra $E_8$.
In fact, that $j$ should encode the representations of $E_8$ was settled \cite{kac} long before the final proof of the Moonshine conjectures \cite{borcherds}.
This relationship between the largest sporadic group and the largest exceptional algebra would connect the McKay Correspondence to Moonshine and thereby weave another beautiful thread into the panoramic tapestry of mathematics.

Over the years, there have been various generalizations of Monstrous Moonshine (cf.~\cite{Gannon:2004xi}) in mathematics and, more recently, in physics \cite{Eguchi:2010ej,He:2012kw,Cheng:2013wca,DGO,Cheng:2015fha}.
Of the vast literature, we will focus on relating the exceptional algebras and the sporadics, review the pertinent background and present a new set of correspondences, which though numeralogically intriguing, remain fundamentally mysterious and await further exploration.

\begin{figure}[th!!!]
\includegraphics[trim=0mm 0mm 0mm 0mm, width=6.2in]{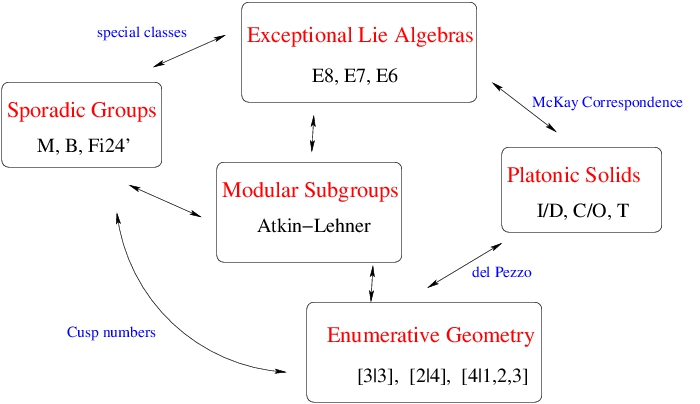}
\caption{
{\sf
The web of correspondences amongst the exceptional Lie algebras and sporadic groups, via some classical enumerative geometry and modular subgroups.
}
\label{f:web}}
\end{figure}

To guide the reader, we summarize the web of inter-connections in Figure \ref{f:web} as well as the ensuing list which point to the relevant sections:
\begin{itemize}
\item McKay Correspondence between $\widehat{E}_{6,7,8}$ Dynkin diagrams and symmetries of Platonic solids in \eqref{affineADE} and \eqref{ADE};
\item 
Dimensions of fundamental representations of $E_{6,7,8}$ in \eqref{dimF}, number of lines (respectively bitangents and tritangents) as well as $(-1)$-curves in classical geometry of del Pezzo surfaces and curves thereon in Table \ref{t:geom};
\item Classes of involutions in $\IM$ and the nodes of the $\widehat{E_{8}}$ Dynkin diagram in \eqref{Me8}, similarly, involutory classes in $2.\IB$ and $\widehat{E_{7}}$ in \eqref{Be7}, and 3-transposition classes of $3.Fi_{24}'$ and $\widehat{E_{6}}$ in \eqref{Fie6};
\item
Cusp number sums of the invariant groups of McKay-Thompson series for $\IM,2.\IB,3.Fi_{24}'$ and the enumerative geometries for the del Pezzo surfaces $dP_{3,2,1}$ in Observations \ref{360}, \ref{448}, and \ref{440};
\item
Conway's group $Co_1$ and the Horrocks-Mumford bundle in Observation \ref{HMCo1}.
\end{itemize}

The outline of the paper is as follows.
We begin with setting the notation for the various subjects upon which we will touch, as well as reviewing the rudiments in as a self-contained manner as possible in \S\ref{s:rev}, before delving into our correspondences in \S\ref{s:cor}. We end with a digression on the Horrocks-Mumford bundle in \S\ref{s:HM}.

\section{Rudiments and Nomenclature}\label{s:rev}\setall
First let us set the notation, and refresh the reader's mind on the characteristics of the various {\it dramatis personae}.

\subsection{$PSL(2,\IZ)$ and $PSL(2,\IR)$}
We will denote the modular group as 
\begin{equation}
\Gamma := PSL(2,\IZ) = \{
{\scriptsize \left(\begin{array}{cc} a & b \\ c & d \end{array} \right)} \  |  \ a,b,c,d \in \IZ, ad-bc=1 \} / \{\pm \II\} \ .
\end{equation} 
The two important subgroups of $\Gamma$ for our purposes are the congruence subgroups:
\begin{align}
\nn
\Gamma_0(N) &:= \{ \gamma \in \Gamma \ | \ c \bmod N = 0 \} \ ;
\\
\Gamma(N) &:= \{ \gamma \in \Gamma \ | {\scriptsize \left(\begin{array}{cc} a & b \\ c & d \end{array} \right)} \bmod N = \II \} \subset \Gamma_0(N) \ .
\end{align}

\paragraph{Cusps: }
The action of $\Gamma$ on the upper-half plane $\cH = \{z : \im(z) > 0\}$ by linear-fractional transformation $z \mapsto \frac{az+b}{cz+d}$ is fundamental.
In addition to $\cH$, of particular interest is the set of {\em cusps} $\IQ \cup \{ \infty \}$ on the boundary real axis.
These are rational points which are taken to themselves under the linear-fractional maps in $\Gamma$ and are the only points on the real axis which should be adjoined to $\cH$ when considering the action of $\Gamma$.
Thus, we will generally speak of the extended upper half plane
\begin{equation}
\cH^* := \cH \cup \IQ \cup \{ \infty \} \ .
\end{equation}
For any subgroup $\Theta$ of $\Gamma$ (including itself), we can define the set of cusps as
\begin{equation} 
C(\Theta) :=  \{ \Theta\mbox{-orbits of } \IQ \cup \{ \infty \} \}.
\end{equation}
An important fact about congruence subgroups $\Theta$ is that $C(\Theta)$ is finite\footnote{
As above, we can show this using representatives modulo the level of the congruence and see the image of $x = \frac{p_x}{q_x}$.
} and
we will henceforth call the number of elements in $C(\Theta)$ the {\bf cusp number}.

Indeed, with respect to the full modular group, for any two points in the extended upper half-plane \footnote{
This can be easily shown by writing, as explicit fractions, $x = \frac{p_x}{q_x}$ and $y = \frac{p_y}{q_y}$ and solving for integer $a,b,c,d$ such that $\frac{p_y}{q_y} = \frac{a \frac{p_x}{q_x} + b}{c \frac{p_x}{q_x} + d}$ and $ad-bc=1$.
} 
$x,y \in \IQ \cup \{ \infty \}$, there exists a $\gamma \in \Gamma$ such that $\gamma(x)= y$. 
Thus, the cusp number of the full modular group is 1, i.e., $|C(\Gamma)| = 1$.

For the Hecke groups, the index of $\Gamma_0(N)$ in $\Gamma$ is $N \prod\limits_{p | N} \left( 1 + \frac{1}{p} \right)$ and
\begin{equation}\label{cuspGamma0N}
|C(\Gamma_0(N)| = \sum\limits_{d|N, d > 0} \phi(\gcd(d,N/d)) \ ,
\end{equation} 
where $\phi(m) = m \prod\limits_{p | m} \left( 1- \frac{1}{p} \right)$ is the standard Euler totient counting the number of positive integers between 1 and $m$, inclusive, which are co-prime to $m$.

\paragraph{Modular Curves: }
By adjoining appropriate cusps which serve as compactification points, we can form the quotient $\Theta \backslash \cH^*$ of the extended upper half plane by any subgroups $\Theta$ of $\Gamma$.
The first classical result \footnote{
The situation can be extended to the reals (cf.~\cite{Gannon:2004xi,Gannong0}).
The group $PSL(2,\IR)$ acts on $\cH$ (here the concept of cusps is irrelevant) and its subgroups $G$ can act similarly. Therefore, one can also form the relevant quotients here.
For example, for $\Gamma$, we have $G \backslash \cH$ being $\IP^1$ with a single puncture (since we did not add the cusp to compactify).
In fact for {\em any} Riemann surface $\Sigma$ of genus $g$ and $n$ punctures with $3g+n > 3$, one can find a subgroup $G \subset PSL(2,\IR)$ such that $G \backslash \cH \simeq \Sigma$ and the fundamental group $\pi_1(\Sigma) \simeq G$.
}
, dating to Klein and Dedekind, is that for the full modular group $\Gamma$, this quotient is the Riemann sphere:
\begin{equation}
\Gamma \backslash \cH^* \simeq \IP^1 \ .
\end{equation}
In general, we can quotient $\cH^*$ by a congruence subgroup and obtain a Riemann surface, dubbed the {\em modular curve}.
Commonly, we denote the modular curves associated to $\Gamma$, $\Gamma_0(N)$ and $\Gamma(N)$ respectively as $X(\Gamma) \simeq \IP^1$, $X_0(N)$ and $X(N)$.

What subgroups $\Theta$ also have $X(\Theta) \simeq \IP^1$?
This so-called {\em genus zero property} is important in several contexts and subgroups possessing it are quite rare.
For example, that there are only 33 (finite-index) torsion-free genus zero subgroups was the classification of \cite{classSebbar} and the physical interpretation, the subject of \cite{He:2012kw,He:2012jn}.

\paragraph{Normalizer: }
Let us consider the normalizer of $\Gamma_0(N)$ in $PSL(2,\IR)$,
\begin{equation}
\Gamma_0(N)^+ := \{
\gamma \in PSL(2,\IR) \ | \ \gamma s = t \gamma\ , \ \exists \ s,t \in \Gamma_0(N) \ .
\}
\end{equation}
We can write the elements of this normalizer more explicitly \cite{CN,AL} as follows.
Let $h$ be the largest divisor \footnote{
The curious fact is that divisors $h$ of 24 are the only integers for which $xy \equiv 1 \bmod h \Rightarrow x \equiv y \bmod h$.
}
of 24 for which $h^2$ divides $N$ and let $N = nh = kh^2$, then
\footnote{
Another way, perhaps more explicit, to write this is
\[
\Gamma_0(N)^+ := \{
\sqrt{e} {\scriptsize
\mat{a & b \\ c & d}} \in PSL(2, \IR)\ , \ | \
a,b,c,d,e \in \IZ, \ ad-bc = e \ , \quad
e | N, \ e | a, \ e|d, \ N | c
\} \ .
\]
}
\begin{equation}
\Gamma_0(N)^+ := \{
\gamma = {\scriptsize
\mat{ae & \frac{b}{h} \\ cn & de} \ , \ | \
\frac{e}{k} \in \IZ \ ,
\gcd(e, \frac{k}{e}) = 1 \ ,
\det(\gamma) = ade^2 - bck = e \in \IZ_{>0}
}
\} \ .
\end{equation}

The above may seem a little difficult to use and present. In order to more conveniently describe it, as will be later needed extensively in the context of Moonshine, we need some further nomenclature.
An important subgroup of $\Gamma_0(N)^+$ is generated by the so-called {\em Fricke involution} $\tmat{0&-1\\N&0}$ taking $z \mapsto -\frac{1}{N z}$.
This extends $\Gamma_0(N)$ to a group 
\footnote{
For primes $p$, using the Fricke involution, we can generate $\Gamma_0(p)^+$ quite simply as:
\[
\Gamma_0(p)^+ = \gen{\Gamma_0(p), \frac{1}{\sqrt{p}}\mat{0&-1\\p&0}} \ .
\]
}
called the Fricke group inside $\Gamma_0(N)^+$, and in which the former is of index 2.
More generally, we have the following
\begin{itemize}
\item
There are the so-called {\bf Atkin-Lehner} involutions \cite{AL} of which Fricke is a special case which are matrices $W_e$ of the form
(here we adhere to the notation of \cite{CN,scheit}) in $\Gamma_0(N)$:
\begin{eqnarray}
W_e = \frac{1}{\sqrt{e}} 
\left( \begin{array}{cc} a & b \\ c & d \end{array} \right) 
\left( \begin{array}{cc} e & 0 \\ 0 & 1 \end{array} \right)
\ , \qquad
\begin{array}{l}
e || N \ , \ \ h = \frac{N}{e} \ ; \\
\scriptsize{\left( \begin{array}{cc} a & b \\ c & d \end{array} \right) }
\in \Gamma_0(h) \ , d \equiv 0 (\bmod e) \ ,
\end{array}
\end{eqnarray}
where the $||$ symbol denotes the {\em Hall divisor}, i.e., $e || N$ means $e | N$ and $\gcd(e, \frac{N}{e}) = 1$.
The set $W_e$ forms a coset of $\Gamma_0(N)$ in $\Gamma_0(N)^+$ and satisfies the relations that, modulo $\Gamma_0(N)$, $W_e^2 = \II$ and $W_eW_f = W_fW_e = W_{\frac{ef}{\gcd(e,f)^2}}$.

\item For $h | n$ and $F_h := {\scriptsize \left(
\begin{array}{cc}
 h & 0 \\
 0 & 1 \\
\end{array}
\right)}$, and define
\begin{align}
\nn
\Gamma_0(n | h) & := F_h^{-1} \Gamma_0(\frac{n}{h}) F_h = 
\left\{
{\scriptsize{
\left(
\begin{array}{cc}
 a & \frac{b}{h} \\
 c h & d \\
\end{array}
\right)} | 
{\scriptsize
\left(
\begin{array}{cc}
 a & b \\
 c & d \\
\end{array}
\right) \in \Gamma_0(\frac{n}{h})}
}
\right\} \ ;
\\
w_e & := F_h^{-1} W_e F_h = \frac{1}{\sqrt{e}} 
\left( \begin{array}{cc} a & \frac{b}{h} \\ ch & d \end{array} \right) 
\left( \begin{array}{cc} e & 0 \\ 0 & 1 \end{array} \right) \ , \quad
m || \frac{n}{h} \ .
\end{align}
The quantities $w_m$ are then the Atkin-Lehner involutions for this group $\Gamma_0(n|h)$.

\item We now introduce the short-hand notation which has become standard to the literature since \cite{CN} (as always, $\left<x_1,x_2\ldots\right>$ denotes the group generated by the elements $x_1,x_2,\ldots$):
{\blue
\begin{equation}
\begin{array}{|l|l|}\hline
n|h+e_1,e_2,\ldots & \left< \Gamma_0(n|h), e_1, e_2, \ldots \right> \\
\hline
n|h+ & \mbox{if all $e || \frac{n}{h}$ are present} \\
\hline
n|h- & \mbox{if all $e || \frac{n}{h}$ are absent}, i.e., n|h-=\Gamma_0(n|h)
\\
\hline
n| &  \mbox{if $h=1$}
\\ 
\hline
\end{array}
\label{n|h+}
\end{equation}
}
\end{itemize}

With the above notation, the key result is the theorem \cite{AL} that
{\purple
\begin{equation}
\Gamma_0(N)^+ = n|h+ \ ; \qquad
N = nh \ , h \mbox{ is the largest divisor of  24~s.t.~} h^2 | N \ .
\end{equation}
}
Indeed, in our shorthand $N|$ simply denotes $\Gamma_0(N)$ and $1$ denotes the full modular group $\Gamma = PSL(2,\IZ)$.

\paragraph{Principal Moduli: }
The central analytic object, again a classical realization dating at least to Klein, is the {\em $j$-invariant}, which is the ``only'' meromorphic function defined on the upper-half plane invariant under the full modular group;
by ``only'' we mean that all invariant functions are rational functions in $j(z)$.
Thus the modular action of $\Gamma$ on $\cH$ leaves invariant the field $\IC(j)$ of rational functions of $j$.
Other than a simple pole at $i \infty$, $j(z)$ is the only holomorphic function invariant under $\Gamma$ once we fix the normalization
\begin{equation}
j(\exp(\frac{2\pi i}{3})) = 0 \ , \quad j(i) = 1728 \ , \qquad
j(\gamma z) = j(z) \ , \gamma \in \Gamma \ .
\end{equation}
Writing the {\bf nome} $q := \exp(2 \pi i z)$, we can perform the famous Fourier expansion of $j(q)$ as
\begin{equation}
j(q) = \frac{1}{q} + 744 + 196884q + 21493760 q^2 + 864299970q^3 + \ldots
\end{equation}
Now, the pole at $z = i \infty$ (i.e., for $q=0$) is explicit and all the coefficients are positive integers.
In the ensuing we will often make use of the {\it normalized} form where the constant 744 has been set to 0, this is habitually denoted as $j_{\IM}$, the subscript will become clear in the following section.
Furthermore, we sometimes divide by 1728 to ensure ramification only at 0,1 and $\infty$.
To clarify our convention, we adhere to the following
{\red
\begin{align}
\nn
j_{arithmetic}(q) & = j(q) = \frac{1}{q} + 744 + 196884q + 21493760 q^2 +\ldots
\\
\nn
j_{analytic}(q) &= \frac{1}{1728} j(q)
\\
j_{\IM}(q) &= j(q)-744 \ .
\end{align}
}

The $j$-invariant is a special case of a {\bf hauptmodul}, or {\it principal modulus}.
For genus zero subgroups $\Theta \subset \Gamma$, the modular functions, i.e., the field of functions invariant under $\Theta$, is generated by a {\em single function}, much like the aforementioned case of the full modular group $\Gamma = PSL(2,\IZ)$ where the $j$-invariant generates $\IC(j)$.
For higher genera, two or more functions are needed to generate the invariants, and there is not as nice a notion of a unique canonical choice \footnote{
Nevertheless, we have such invariants due to Igusa and Shioda and the reader is referred to \cite{Bose:2014lea} for realizations thereof in physics as well as the references therein.
}.
There is, however, at least a notion of replicable functions for higher genus and the reader can consult the nice work of Smith \cite{smith}.

In general, a principal modulus for an arbitrary subgroup $\Theta$ of $\Gamma$ can be seen as an isomorphism from $\Theta \backslash \cH$ to $\IC$, normalized so that its Fourier series starts as $q^{-1} + \cO(1)$.
Now, for any genus, up to conjugation there are only a finite number \cite{thompson} of subgroups of $PSL(2,\IR)$ and for our case of genus 0, there are \cite{cummins} precisely 6484.
Of these, 616 infinite series have integer q-expansion coefficients (cf.~\cite{gGT1} for some recent work beyond genus 0 as well as a classic work on why genus 0 in \cite{CCS}).

\subsection{The Monster}
Of the many fascinating properties of the Monster sporadic group, we will make particular use of the following, in the context of the {\bf Moonshine Conjectures} \cite{CN}, some initial computations of which were settled in \cite{AFS,FLM} and much of which later proven in \cite{borcherds}.

\paragraph{Supersingular Primes: }
The order of the Monster is
{\purple
\begin{equation}\label{|M|}
|\IM| =  2^{46} \cdot 3^{20} \cdot 5^9 \cdot 7^6 \cdot 11^2 \cdot 13^3 \cdot 17 \cdot 19 \cdot 23 \cdot 29 \cdot 31 \cdot 41 \cdot 47 \cdot 59 \cdot 71 
\sim 10^{54}
\ .
\end{equation}
}
The observation of Ogg \cite{Ogg} was that a prime $p$ such that $\Gamma_0(p)^+ \backslash \cH$ is genus zero {\em if and only if} $p$ is one appearing in the above list.
To this day, though the Moonshine conjectures \cite{CN} have been proven \cite{borcherds}, this earliest observation on Moonshine remains unexplained and Borcherds' proof does not actually contain a direct explanation of the appearance of these particular primes.
Thus, Ogg's prize of a bottle of Jack Daniels is yet to be collected \cite{DO}.

These primes \footnote{
Incidentally, all these primes are Chen primes, i.e., primes $p$ such that $p+2$ is either itself a prime, or the product of exactly 2 primes.
} are called {\it supersingular} and in summary, they are the 15 primes obeying the following equivalent definitions:
\begin{itemize}
\item The modular curve $X_0(p)^+ = \Gamma_0(p)^+ \backslash \cH$ is genus zero (where $\Gamma_0(p)^+$ is the congruence group adjoining the Fricke involution as defined earlier);
\item 
The terminology ``supersingular'' coincides with that in the theory of elliptic curves for a reason: every supersingular elliptic curve over $\IF_{p^r}$ can be in fact defined just over the subfield $\IF_p$.
The simplest working definition of such a curve is that, in Legendre form
$y^2 = x(x-1)(x-\lambda)$, we have that (cf.~\cite{BM}) the Hasse invariant
$\sum\limits_{i=0}^{\frac12(p-1)} {\frac12(p-1) \choose i}^2 \lambda^i = 0$;
\item 
In \cite{pizer}, it was noticed that the Hecke Conjecture is true for $\Gamma_0(p)$ only for these primes
\footnote{
The conjecture, proved by Pizer, is that a certain family of theta-series (cf.~\cite{pizer}) for prime $p$ and a quaternionic division algebra is linearly independent only for $p$ supersingular.
}.
\item 
The second author later observed \footnote{
In a recent work, Erdenberger \cite{Erden} points out another extraordinary emergence of these very same primes as follows.
Consider the Siegel upper plane $\IH_2 := \left\{
\tau = \left(\begin{array}{cc}
\tau_1 & \tau_2 \\ \tau_2 & \tau_3 \\
\end{array}\right) \in \mbox{Sym}^2(2;\IC),\ \im(\tau) > 0
\right\}$, and a subgroup $\Gamma_p$ of the symplectic group as
$\Gamma_p := \left\{
M \in Sp(4;\IQ) : M \in {\tiny \left(\begin{array}{cccc}
\IZ & \IZ & \IZ & p \IZ \\
p \IZ & \IZ & p \IZ & p \IZ \\
\IZ & \IZ & \IZ & p \IZ \\
\IZ & \frac{1}{p} \IZ & \IZ & \IZ \\
\end{array}\right)}
\right\}$.
The action of $\Gamma_p$ on $\IH_2$ can be given as $M = \left(
\begin{array}{cc} A & B \\ C & D
\end{array}\right) : \tau \mapsto (A\tau + B)(C\tau+D)^{-1}$ with $A,B,C,D$ two by two matrices acting on the matrix $\tau$.
Thus defined, a non-trivial cusp (modular) form for $\Gamma_p$ exists if $p > 71$ or if $p \in \{37,43,53,61,67\}$.
The complement of these is precisely the set of primes which appear in \eqref{|M|}, dividing the order of the Monster, the same as the Ogg list; a clarification of this is under way \cite{ogglist}.
}
that only trivial cusp-form for the symplectic group $\Gamma_p$ as classified in \cite{Erden} exists when $p$ is one of the above 15 primes;
\item
Recently, Duncan and Ono \cite{DO} point out that the McKay-Thompson series (which we will introduce in detail in \eqref{MT}) are encoded by precisely the $j$-invariants of supersingular elliptic curves. 
\end{itemize}

It is interesting that for $\Gamma_0(N)$ itself, one has $X_0(N) = \Gamma_0(N) \backslash \cH$ being genus zero precisely for 15 values of $N$, namely
\begin{equation}\label{g=0N}
N = 1, 2, 3, 4, 5, 6, 7, 8, 9, 10, 12, 13, 16, 18, 25 \ .
\end{equation}
These values have recently \cite{Cheng:2013wca} been pointed out to match the Coxeter numbers of root systems of Niemeier Lattices (cf.~Eq 2.20 therein).

\paragraph{Moonshine: }
There are 194 (linear, ordinary) irreducible representations of $\IM$, and thus 194 conjugacy classes.
This constitutes a standard $194 \times 194$ character table, the first column of which is the vector of the dimensions of the irreducible representations, starting with
\begin{equation}
\{
1, 196883, 21296876, 842609326, 18538750076, \ldots
\} \ .
\end{equation}
These are the dimension (degree) of the irreducible representations $\rho_1$, $\rho_{196883}$, $\rho_{21296876}$, etc.

Now, consider an infinite-dimensional representation of $\IM$
\begin{equation}
V = V_0 \oplus V_1 \oplus V_2 \oplus \ldots
\end{equation}
with $V_0 = \rho_1$, $V_1 = \{0\}$, $V_2 = \rho_1 \oplus \rho_{196883}$,
$V_3 = \rho_1 \oplus \rho_{196883} \oplus \rho_{21296876}$, \ldots
The corresponding generating function (graded dimension) is then
\begin{equation}
\sum\limits_{n=0}^\infty q^n \dim(V_n) = 
1 + 196884q^2 + 21493760 q^3 = q (j(q) - 744) = q j_{\IM}(q)\ .
\end{equation}
The second equality is remarkable and is part of the key results of Moonshine, relating finite groups to modular groups.

Now, a character of an element, $g$, is rational if $g$ is conjugate to
its inverse. Of the 194 conjugacy classes of $\IM$, there are 22 characters which are complex quadratic valued.
If we replace a complex irreducible representation, $R$, by its sum $R + \overline{R}$ and remove duplicate representations and classes, this yields the
rational character tables of $\IM$, on which we now focus.
The vector of dimensions of irreducible representations is the first column of the character table which yields the above sum to the $j$-function.
We can perform a similar sum for all the 172 rational conjugacy classes and obtain a generating function for each as
{\red
\begin{align}\label{MT}
\nn
T_g(q) &= q^{-1} \sum\limits_{n=1}^\infty \ch_{V_n}(g) q^n \\
&=  q^{-1} + 0 + h_1(g)q + h_2(g) q^2 + \ldots
\ ,
\end{align}
}
where $\ch$ are the characters (indeed, characters, being traces of finite matrices, are defined over conjugacy classes) of these representations, called {\it head characters} $H_n(g)$ and $h_n(g) = \tr H_n(g)$. 
This is the {\bf McKay-Thompson series}.
Indeed, for $g = \II$, $T_\II(q)$ is the above (normalized) $j$-function, or $j_{\IM}$ in our notation.
\comment{
Of course, we desire integer coefficients.
The character table of $\IM$ consists almost entirely of integers but there are exceptions.
The non-integer entries to the character table are quadratic and complex; there are 22 two-by-two blocks of these.
Therefore, there are $194-22=172$ integer columns.
}
As we will see in \S\ref{s:cor}, two conjugacy classes (27A and 27B) give rise to the same McKay-Thompson series, hence only 171 (including the identity) are candidates for Moonshine.

There are further numerological mysteries \cite{novel} surrounding $\IM$.
In addition to the supersingular primes mentioned above, the fact that of the 171 integer characters, some of the associated McKay-Thompson series are linearly dependent over $\IZ$, and that subsequently there are 163 $\IZ$-independent McKay-Thompson series for the Monster (cf.~pp 310 and 317 of \cite{CN}).
This is intriguing:
it is well known that 163 is the largest of the {\bf Heegner} numbers.
We recall that there are precisely 9 of such numbers:
\begin{equation}
\mbox{Heegner} = \{1,2,3,7,11,19,43,67,163\} \ .
\end{equation}
These numbers $H$ have the distinction that the imaginary quadratic field $\IQ(\sqrt{-H})$ has class number one, i.e., the associated ring of integers has unique factorization.
\comment{
A nice historical anecdote is the fact that before the Monster was fully constructed, it was thought that it had precisely 163 conjugacy classes.
}
A more recent observation of the second author (cf.~an account in \cite{DW}) is in $E_8 \times E_8$ heterotic string theory compactification on a K3 surface which is dual to F-theory  compactification on a Calabi-Yau threefold elliptically fibred over a complex surface $B$. In the extremal case where one of the $E_8$ gauge groups is completely broken, the base surface $B$ has Picard number exactly 194.
These and ever-increasing number of observations continue to intrigue us \cite{GN,LYY,Duncan,HLY,Gannon:2004xi,DW,Cheng:2015fha,DGO,Cheng:2013wca,He:2014uma}.

\subsubsection{Monstrous Moonshine}\label{s:MM}
The key result of Moonshine is that the McKay-Thompson series defined in \eqref{MT}, for each conjugacy class of an element $g$ in the Monster, has the following property:
\begin{theorem} [Moonshine]
The $q$-series $T_g(q)$ is the normalized generator of a genus zero function field arising from a group between $\Gamma_0(N)$ and its normalizer $\Gamma_0(N)^+$ in $PSL(2,\IR)$. 
\end{theorem}

The integer $N$ can be determined in several equivalent ways \cite{CN}.
Let $\cF(g)$ be precisely the elements of $PSL(2,\IR)$ which fix $T_g$.
That is, as $J(q)$ is {\it the} modular invariant of $\Gamma$, $T_g(q)$ is {\it the} invariant of $\cF(g)$), then 
\begin{itemize}
\item 
$N$ is the {\sf level} of $\cF(g)$;
\item 
$N$ is the {\it smallest} integer so that the group element sending $z \mapsto \frac{z}{Nz+1}$ is in $\cF(g)$;
\item
If $n$ is the order of the (conjugacy class of) group element $g$ in $\IM$, then
\begin{equation}
N / n = h \in \IZ_{>0} ; \qquad h | 24 \ , \quad h^2 | N \ .
\end{equation}
\end{itemize}
The ``Euler characteristic'' (to be detailed later) and cusp number of all the $\cF(g)$ for the 194 conjugacy classes of the Monster were calculated in \cite{CN}.
More recently, a generalization was performed where the following set 
$\Delta := \{ G : genus(G) = 0, \ \Gamma_0(m) \subseteq G \subseteq \Gamma_0(m)^+ \}$ for some integer $m$ is analysed \cite{CL,CMS}.
Here $G$ is some modular subgroup residing between the congruence subgroup $\Gamma_0(m)$ and its normalizer $\Gamma_0(m)^+$ in $PSL(2,\IR)$ and $genus(G)$ is the genus of the Riemann surface $G \backslash \cH$.
Setting $m=nh^2$ where $h$ is the largest divisor of 24 such that $h^2|m$, the number of distinct pairs $(n,h)$ is 419.

\comment{
\subsubsection{Summary: Moonshine Groups and Principal Moduli}
It is expedient to summarize the key properties of the functions and groups arising from Moonshine, emphasizing the relevant computational details:
\begin{itemize}
\item $G$ is a finite index, genus zero group (i.e., $(G \backslash \cH ) \cup \{\mbox{cusps}\} \simeq \IP^1$) of $SL(2,\IR)$
  \begin{itemize}
  \item Commensurability: $G \cap SL(2,\IZ)$ is finite index in both $G$ and $SL(2,\IZ)$;
  \item Contains congruence group: $\Gamma_0(N) \subset G \subset SL(2,\IR)$;
  \item meromorphy at cusps: $\left(\begin{array}{cc} 1 & t \\ 0 & 1 \end{array}\right) \in G$ iff $t \in \IZ$;
  \end{itemize}
\item For conjugacy class $\cC$ on $\IM$, 
\end{itemize}

}

\subsection{Exceptional Affine Lie Algebras}
The Cartan-Killing classification of simple Lie algebras is a triumph of late C19th mathematics.
Of the Dynkin diagrams, the simply-laced ones consist only of single bonds and fall into an ADE pattern: the two infinite series $A_n \simeq \mathfrak{sl}_{n+1}(\IC)$ and $D_n \simeq \mathfrak{so}_{2n}(\IC)$, and the exceptionals $E_{6,7,8}$, the Dynkin diagrams of which are:
\[
\begin{array}{c}
\begin{array}{l}\includegraphics[trim=0mm 0mm 0mm 0mm, width=5.0in]{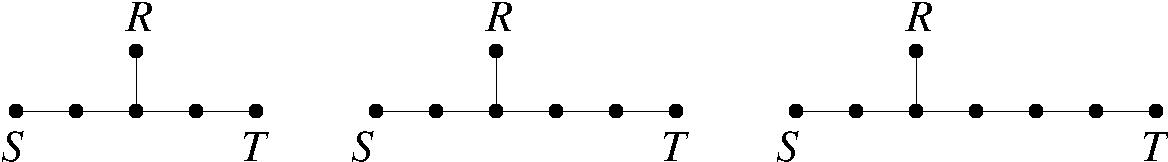}\end{array}\\
E_6 \hspace{2in} E_7 \hspace{2in} E_8
\end{array}
\]
To each of the above is associated a Platonic solid(s), where the famous five group themselves exactly into 3 classes in the sense that the cube and the octahedron, as well as the dodecahedron and the icosahedron are graph duals and share the same symmetry group, whereas the tetrahedron is self-dual:
\[
\begin{array}{|c|c|c|}\hline
E_6 \sim \mbox{Tetrahedron} &
E_7 \sim \mbox{Cube/Octahedron}&
E_8 \sim \mbox{Dodecahedron/Icosahedron} 
\\
\includegraphics[trim=0mm 0mm 0mm 0mm, height=1in]{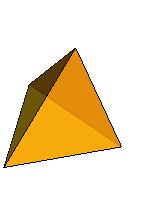}
&
\includegraphics[trim=0mm 0mm 0mm 0mm, height=0.6in]{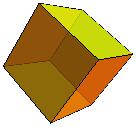}
\includegraphics[trim=0mm 0mm 0mm 0mm, height=1in]{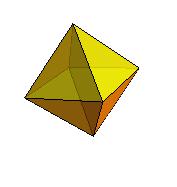}
&
\includegraphics[trim=0mm 0mm 0mm 0mm, height=1in]{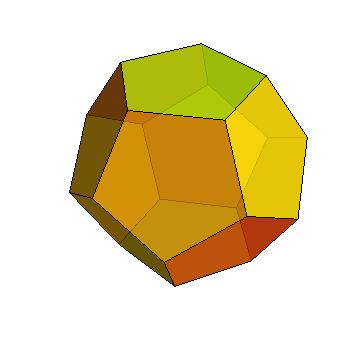}
\includegraphics[trim=0mm 0mm 0mm 0mm, height=1in]{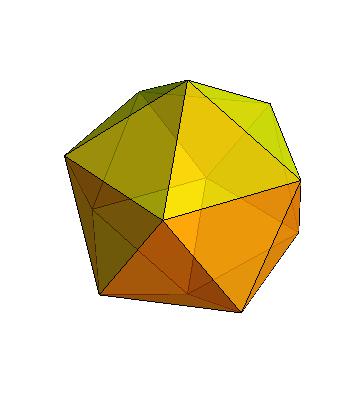}
\\
\hline
\end{array}
\]
The finite group of symmetries of each solid can be read off from the Dynkin diagram as follows.
Associate generators $R,S,T$ to the extremalities of the Dynkin diagram and the order thereof is equal to the number of bonds between it and the central trivalent node, plus 1:
\begin{equation}\begin{array}{rcl}\label{adegroups}
E_6 & : & G = \gen{R,S,T \ | \ RST=R^2 = S^3 = T^3 = \II} \simeq A_4 \ , \quad |G| = 12 \ , \\
E_7 & : & G = \gen{R,S,T \ | \ RST=R^2 = S^3 = T^4 = \II} \simeq S_4 \ , \quad  |G| = 24 \ , \\
E_8 & : & G = \gen{R,S,T \ | \ RST=R^2 = S^3 = T^5 = \II} \simeq A_5 \ , \quad  |G| = 60 \ . \\
\end{array}\end{equation}
Incidentally, We see that in each case $|G|$ divided by the order of $R,S,T$ gives a triple, viz., $(6,4,4)$, $(12,8,6)$ and $(30,20,12)$. We recognize these as the number of (edges, vertices/faces, faces/vertices) of the corresponding solids.

The relationship 
\footnote{
There is another intriguing observation of Kostant \cite{kostant}.
The group $PSL(n,q)$ of uni-determinant $n\times n$ matrices over the finite field of $q$ elements is a finite group, it acts non-trivially on the projective space $\IP^n(\IF_q)$ which has $\frac{q^n-1}{q-1}$ elements.
However, it is rare that it acts only on a strict subset of these elements.
For $q=p$, some prime, this only happens for $PSL(n,p)$ at $p=2,3,5,7,11$ and the group acts non-trivially only on $p$ points.
Of these five, only when $p=5,7,11$ is $PSL(n,q)$ a simple finite group and in fact does not act non-trivially on fewer than $p$ points (a fact known to Galois).
Remarkably, 
$PSL(2,5) \simeq A_4 \times_{set} \IZ_5$, 
$PSL(2,7) \simeq S_4 \times_{set} \IZ_7$, and
$PSL(2,11) \simeq A_5 \times_{set} \IZ_{11}$,
where the notation $\times_{set}$ is to emphasize that it is not a group product (after all, these groups are simple).
We see the emergence of $E_{6,7,8}$ here from the factors $A_4$, $S_4$ and $A_5$.
} between the Lie groups and the Platonic solids was made more striking in \cite{mckay}, where by taking the double covers of the groups in \eqref{adegroups} - thus making them discrete subgroups of $SU(2)$, rather than $SO(3)$ in which the rotational symmetries of the solids are visualized - and distinguishing the fundamental 2-dimensional complex irreducible representation $\cR$.
These {\it binary} groups \cite{CCS} will have orders twice those in \eqref{adegroups}:
{\purple
\begin{equation}\label{affineADE}\begin{array}{c|c|c}
\mbox{Group} & \mbox{Presentation} & \mbox{Order} \\ \hline
\widehat{E_6} & \{r,s,t \ | \ r^2 = s^3 = t^3 = rst \} & 24 \\
\widehat{E_7} & \{r,s,t \ | \ r^2 = s^3 = t^4 = rst \} & 48 \\
\widehat{E_8} & \{r,s,t \ | \ r^2 = s^3 = t^5 = rst \} & 120
\end{array}\end{equation}
}
Note that the only difference in the presentation is that we remove the condition $=\II$ in the definitions.
Subsequently, the operator $\cR \otimes$ gives the decomposition over the irreducible representations $\{ R_i \}$:
\begin{equation}\label{ADE}
\cR \otimes R_i = \bigoplus\limits_j a_{ij} R_j \ .
\end{equation}
Remarkably, the $a_{ij}$ matrices are precisely the adjacency matrices \footnote{
These adjacency matrices are also exactly those with maximum eigenvalue two \cite{smithJ}.
} of the {\it affine} Dynkin diagrams of the {\it extended} or {\it affine} semi-simple Lie algebra of $\widehat{ADE}$ type (the gauge theory implication of this is discussed in \cite{Hanany:1998sd}); this is the {\it McKay Correspondence}.

In this correspondence, each node of the affine Dynkin diagram is associated to an irreducible representation of the corresponding group: the dimension of the irreducible representation precisely matches the dual Coxeter labels, which are the expansion coefficients $a_i^\vee$ of the normalized highest root $\theta$ into the basis $\{\alpha^{(i)\vee}\}$  of simple coroots: $\frac{2}{(\theta,\theta)} \theta = \sum\limits_{i=1}^r a_i^\vee \alpha^{(i)\vee}$.
Being dimensions of irreducible representations, the sum of squares of these labels will be precisely the orders of the associated binary groups: 24, 48 and 120.
The Dynkin diagrams are (with the affine node in white)
\begin{equation}\label{affineDynkin}\begin{array}{c}
\begin{array}{l}\includegraphics[width=6.0in]{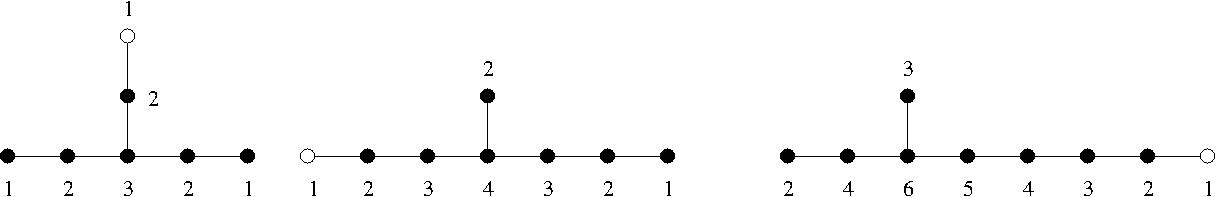}\end{array}\\
\widehat{E_6} \hspace{2in} \widehat{E_7} \hspace{2in} \widehat{E_8}
\end{array}
\end{equation}

Algebro-geometrically, the binary discrete subgroups of $SU(2)$ in \eqref{affineADE} furnish affine models for K3 surfaces as orbifolds of the form $\IC^2 / \Gamma$ and are called du Val singularities \cite{duval}.
They can be described as affine equations in $\IC[x,y,z]$ as follows:
\begin{equation}
\begin{array}{c|c|c}
\mbox{du Val} & \mbox{Defining Eq} & \mbox{Deg}(x,y,z) \\ \hline
E_6 & x^2+y^3+z^4=0 & (6,4,3) \\ \hline
E_7 & x^2+y^3+yz^3=0 & (12,8,6) \\ \hline
E_8 & x^2+y^3+z^5=0 & (30,20,12) \\
\end{array}
\end{equation}
In the above, Deg$(x,y,z)$ means a weight which we can assign to the variables $(x,y,z)$ respectively so that the equations become homogeneous, of degree respectively 12,24,60, the orders of the finite groups $E_{6,7,8}$ for the Platonic solids. 
Comparing with \eqref{adegroups} we could see the remnants of the $R,S,T$ generators and their relations with the number of (edges, faces/vertices, vertices/faces) for each of the solids \footnote{
To get from the solids to these polynomials in a quick way, q.v.~Baez's short introduction to ADE theory in \url{http://math.ucr.edu/home/baez/ADE.html}.
}.

Finally, the numbers relevant to us come from the Lie algebras themselves \footnote{
The reader is referred also to the so-called {\em Arnol'd Trinities} \cite{arnold}, a mysterious web of correspondences involving, inter alia, $E_{6,7,8}$ and $\IR, \IC, \IH$; therein is nice recasting of $24, 48, 120$ in terms of the real projective plane.
Moreover, the numbers 120 and $2\cdot248$ emerge in the context of 2-local subgroups of $\IM$ and $\IB$ \cite{meier}.
}.
We recall that the dimensions of the {\em fundamental representations} of our (ordinary non-affine) exceptional algebra are
\begin{align}
\nn
\dim_F(E_6) &= 27 \ , \\
\nn
\dim_F(E_7) &= 56 = 28 \times 2 \ , \\
\dim_F(E_8) &= 248 = 120 + 120 + 8  \ .
\label{dimF}
\end{align}

\newpage

\subsection{Classical Enumerative Geometry}\label{s:enum}
Parallel to the aforementioned Lie algebras, another set of celebrated C19th mathematics comes from enumerative geometry; these will also be of concern to us (see Hitchin's lectures \cite{hitchin}).
In particular, three counting problems distinguish themselves to us.
We will use the notation that $[n|a_1,a_2,\ldots,a_k]$ means the (not necessarily complete) intersection of $k$ polynomials of degrees $a_1, \ldots, a_k$ respectively in $\IP^n$.
\begin{description}
\item[Cayley-Salmon (1849): ] 
The cubic surface in $\IP^3$, $[3|3]$,  has exactly 27 lines. 

\item[Jacobi (1850): ]
The quartic curve in $\IP^2$ has exactly 28 bitangents.
We recall that bitangents are lines tangent to a curve at 2 different points; indeed, by degree count using B\'ezout, starting at degree 4, curves can have such bitangents. Of course, there are always infinite number of secants and tangents.
Such a curve can be realized as $[2 |4]$ and is a Riemann surface of genus $\frac12(4-1)(4-2) = 3$ by adjunction \cite{hartshorne}.
 
\item[Clebsch (1863): ]
The canonical sextic curve of genus 4 has exactly 120 tritangent planes (i.e., planes which are tangent to the curve at precisely 3 points).
This curve can be realized as $[4 |1,2,3]$, i.e., the intersection of a line, a quadric and a cubic in Fermat form in the 5 homogeneous coordinates of $\IP^4$, giving us the so-called {\it Bring's curve} \cite{edge}.
\end{description}
Specifically, Bring's curves can be realized as the Fermat cubic, sliced by the Fermat quadric, and then the line, in the homogeneous coordinates of $\IP^4$:
\begin{equation}\label{bring}
{\cal B} = \{
\sum_i x_i^3 = \sum_i x_i^2 = \sum_i x_i = 0
\} \subset \IP^4 \ .
\end{equation}

These classic results may at first seem esoteric.
However, our attention is drawn to the numbers 27, $56 = 28 \cdot 2$ and $360 = 120 \cdot 3$.
This is not a coincidence and is well-understood in terms of del Pezzo surfaces (cf.~\cite{dolgachev}).
Now, it is well known that the second homology $H_2(dP_d;\IZ)$ of a del Pezzo surface of degree $d$ is generated by the hyperplane class $H$ on the $\IP^2$, as well as the $9-d$ exceptional blow-up $\IP^1$-curve classes.
The intersection matrix of these $d+1$ classes is the Cartan matrix of the affine $\widehat{E_{9-d}}$ algebra and whence the adjacency matrix of the associated Dynkin diagram (cf.~\eqref{ADE}).

Indeed, $[3|3]$ is birational to $\IP^2$ blown up at 6 generic points \cite{hartshorne}, furnishing a del Pezzo surface $dP_3$ of degree 3, and there are 27 lines pass through these blow-up points appropriately as $(-1)$-curves (curves whose self-intersection equals $-1$).
Likewise, $dP_2$, the del Pezzo surface of degree 2, has 56 $(-1)-$curves.
The linear system of the anti-canonical divisor of $dP_2$ maps to $\IP^2$ branched over $[2 | 4]$ and the 56 curves pair to the 28 bitangents.
Finally, $dP_3$, has 240 $(-1)$-curves. The linear system of its anti-canonical divisor maps to $\IP^2$ branched over $[4 |1,2,3]$, with these 240 curves pairing to the 120 tritangents (NB.~the order of the binary icosahedron group is 120).
Furthermore, we can explicitly see the Weyl groups of the root system of the respective Lie algebras are the automorphism groups of the aforementioned geometric objects \cite{hitchin,manivel}.

Of interest to us also, since we are touching on the subject of bitangents, is the {\bf theta characteristic} of an algebraic curve $X$.
We recall \cite{coble,dolgachev} that this is an element $\vartheta \in Pic(X)$, the Picard group of line bundles on $X$, which squares to the canonical bundle:
 $\vartheta^{\otimes 2} = \omega_X$. It is even/odd according to whether the number of global sections $h^0(X,\vartheta)$ is even/odd.
We have that the number of theta charactertistics on $X$ of genus $g$ is
\begin{equation}\label{theta}
\#(\mbox{odd } \vartheta) = 2^{g-1}(2^g-1) \ , \quad
\#(\mbox{even } \vartheta) = 2^{g-1}(2^g+1) \ .
\end{equation}
The total number $2^{2g}$ is, incidentally, the number of points of the Jacobian $Jac(X)$ defined over the finite field $\IF_2$.
Importantly, the number of bitangent planes to a curve $X$ is precisely that of odd theta characteristics in \eqref{theta}.

In summary, for the above three classical enumerative problems and in relation to the del Pezzo surfaces, we collect the relevant facts in Table \ref{t:geom}.

\begin{table}[h!!!]
\[
\begin{array}{|c|c|c|c|c|c||c|c|}\hline
\mbox{Geometry} & \mbox{Configuration } \cC& \mathfrak{g} & S = W(\mathfrak{g}) = Aut(\cC) & |S| & dP_d & (-1)\mbox{-Curves} \\
\hline \hline
[3|3] & \mbox{27 Lines} & \widehat{E_6} & Aut(PSU_4(2)) & 2^7 \cdot 3^4 \cdot 5 
& dP_3 & 27 
\\ \hline
[2|4] & \begin{array}{c} \mbox{28 Bitangent} \\ \mbox{Lines} \end{array} 
      & \widehat{E_7} & \IZ_2 \times PSp_6(2) & 2^{10} \cdot 3^4 \cdot 5 \cdot 7
& dP_2 & 56 
\\ \hline
[4|1,2,3] & \begin{array}{c} \mbox{120 Tritangent} \\ \mbox{Planes} \end{array}
      & \widehat{E_8} & O_8(2)^+ & 2^{14} \cdot 3^5 \cdot 5^2 \cdot 7
& dP_1 & 240 
\\ \hline
\end{array}
\]
{\caption
{\sf
The correspondences between 3 classical enumerative geometrical problems and the exceptional Lie algebras.
The geometry, $[n|a_1,a_2,\ldots,a_k]$ means the intersection of $k$ polynomials of degrees $a_1, \ldots, a_k$ respectively in $\IP^n$. $W(\mathfrak{g})$ means the Weyl group of the root system of the Lie algebra $\mathfrak{g}$, which here is equal to the automorphism group $Aut(\cC)$ of the configuration of lines and tangents in the geometry.
}
\label{t:geom}
}
\end{table}

\newpage

\section{Correspondences}\label{s:cor}\setall
Let us first examine the 194 conjugacy classes of $\IM$ in more detail.
In the standard notation of ATLAS \cite{atlas}, the classes are recorded as `$n$X' where $n$ is the order of an element in the class and X is a capital letter indexing the classes of order $n$, ascending alphabetically according to increasing sizes of the centralizer.

In particular, the Monster has classes 1A; 2A, 2B; 3A, 3B, 3C, \ldots, 110A; 119A, 119B.
There are 73 distinct orders, and below we write the multiplicities as subscripts (over which the letters run)
\begin{align}
\nn
\mbox{Classes}(\IM) 
= \{
&
1_1,2_2,3_3,4_4,5_2,6_6,7_2,8_6,9_2,10_5,11_1,12_{10},13_2,14_3,15_4,16_3,17_1,18_5,\\ \nn
&
19_1,20_6,21_4,22_2,23_2,24_{10},25_1,26_2,27_2,28_4,29_1,30_7,31_2,32_2,33_2,
\\ \nn
& 34_1,35_2,36_4,38_1,39_4,40_4,41_1,42_4,
   44_2,45_1,46_4,47_2,48_1,50_1,51_1,52_2,54_1,
\\ \nn
& 55_1,56_3,57_1,59_2,60_6,62_2,66_2,68_1,69_2,70_2,71_2,78_
   3,84_3,87_2,88_2,92_2,
\\ \nn
& 93_2,94_2,95_2,104_2,105_1,110_1,119_2
\} \ .
\end{align}

Of the 194 classes, each giving a column in the standard character table, 22 are obviously grouped as complex conjugates (i.e., Galois orbits)\footnote{
In GAP \cite{gap}, these can readily be found using the {\sf ClassOrbit(~)} command for {\sf CharacterTable(``M'')}.
}, these are
\begin{equation}
\begin{array}{c}
(23A,23B); \
(31A,31B); \
(39C,39D); \
(40C,40D); \
(44A,44B); \
(46A,46B); \
(46C,46D); \\
(47A,47B); \
(56B,56C); \
(59A,59B); \
(62A,62B); \
(69A,69B); \
(71A,71B); \
(78B,78C); \\
(87A,87B); \
(88A,88B); \
(92A,92B); \
(93A,93B); \
(94A,94B); \
(95A,95B); \\
(104A,104B); \
(109A,109B); \
\end{array}
\end{equation}
As mentioned earlier, this gives the rational character table of size $194-22=172$.
With 9 further linear relations amongst the McKay-Thompson series (cf.~p310 of \cite{CN}), we have the column rank of $172-9=163$, the largest Heegner number.

\subsection{Desire for Adjacency}
The second author's initial observation was that the $j$-function not only encodes the irreducible representations of $\IM$ but also that
\begin{equation}
j(q)^{\frac13} = q^{-\frac13}\left(
1 + 248 q + 4124 q^2 + 34752 q^3 + \ldots
\right)
\end{equation}
encodes the irreducible representations of the Lie algebra $E_8$ in a similar fashion (note that $248 = 744/3$):
\begin{equation}
248 = 248 \ , \quad
4124 = 3875 + 248 + 1 \ , \quad
34752 = 30380 + 3875 + 2 \cdot 248 + 1 \ldots
\end{equation}
This was in fact the first puzzle to be settled \cite{kac}:
the unique level-1 highest-weight representation of the affine Kac-Moody algebra $E_8^{{(1)}}$  has graded dimension encoded by $j(q)^{\frac13}$.
One should also be mindful 
\footnote{
Incidentally, the reader is also alerted to the curiosity that $\sigma_1(240) = 744$.
}
of the fact that the theta-series for the $E_8$ root lattice $\Lambda(E_8)$ is $\theta_{\Lambda(E_8)}(q) = \sum\limits_{x \in \Lambda(E_8)} q^{|x|^2/2} = 1 + 240 \sum\limits_{n=1}^\infty \sigma_3(n) q^{2n} = E_4(q)$, the 4th Eisenstein series, so that we have
\begin{equation}\label{thetaE8}
j(q) = \frac{\theta_{\Lambda(E_8)}(q)^3}{\Delta(q)} \ , \qquad
\Delta(q) = \eta(q)^{24} \ , 
\end{equation}
where $\eta(q)$ is the Dedekind eta-function and $\Delta(q)$ is the Ramanujan delta-function.

\subsubsection{Initial Observation on $\IM$ and $\widehat{E_8}$}
Many further curious connections between $\IM$ and $E_8$ have been observed, which still eludes us today (cf.~\cite{GN,LYY,Duncan,HLY,Gannon:2004xi,DW,lam}); to this we will add another here.
The now classic one is as follows.
Consider the two order-2 conjugacy classes $2A$ and $2B$.
The first has $2^4\cdot3^7\cdot5^3\cdot7^4\cdot11\cdot13^2\cdot29\cdot41\cdot59\cdot71 \sim 10^{20}$ elements while the second exceeds it by 7 orders of magnitude.
If we were to multiply any two elements of $2A$, the resulting element can only be in one of 9 conjugacy classes, viz.,
$1A, 2A, 3A, 4A, 5A, 6A, 4B, 2B, 3C$.
The second author then noticed that we have seen these 9 numbers before \cite{FM,GN}!
Glancing back at $\widehat{E_8}$ in \eqref{affineDynkin}, we see that they are precisely the (dual Coxeter) labels of the 9 nodes in the Dynkin diagram, which we also know to be the dimensions of the irreducible representations of the binary icosahedral group by \cite{mckay}.
That is, we have
\begin{equation}\label{Me8}
\begin{array}{ccccccccccccccc}
&&&&&&&& &&3C &&&& \\
&&&&&&&& &&| &&&& \\
1A & - & 2A & - & 3A & - & 4A & - & 5A & - & 6A & - & 4B & - & 2B \\
\end{array}
\end{equation}
The edges, i.e., the meaning of adjacency, however, have no clear interpretation in this correspondence and still awaits clarification \cite{GN,LYY,Duncan,HLY}.
A recent work \cite{Duncan} nicely recasts this observation solely in terms of the properties of $PSL(2,\IR)$.

\subsubsection{The Baby and $\widehat{E_7}$}
\label{s:e7}
An important subgroup of the Monster is the affectionately named Baby Monster, $\IB$, of order $2^{41}\cdot 3^{13}\cdot 5^6\cdot 7^2\cdot 11\cdot 13\cdot 17\cdot 19\cdot 23\cdot 31\cdot 47$.
Its double cover, $2.\IB$, is the centralizer of class $2A$ in $\IM$.
Indeed, comparing with Table~\ref{t:cent}, we see that the order of $2.\IB$ is the size of the associated centralizer of class $2A$.
The observation in \eqref{Me8} was generalized by \cite{HLY} to relate $\IB$ to $E_7$ using the explicit embedding of the vertex algebra.
In summary, we have the product of two involution classes of the Baby falling into 8 classes whose orders coincide with the dual Coxeter numbers of affine $E_7$:
\begin{equation}\label{Be7}
\begin{array}{ccccccccccccc}
&&&& &&2c &&&&&& \\
&&&& &&| &&&&&& \\
1a & - & 2b & - & 3a & - & 4b & - & 3a & - & 2b & - & 1a \\
\end{array}
\end{equation}

\subsubsection{Fischer and $\widehat{E_6}$}\label{s:e6}
Another important subgroup of $\IM$ is the largest of the Fischer groups, $Fi_{24}'$ (sometimes denoted as $F_{3+}$) , of order $2^{21}\cdot3^{16}\cdot5^2\cdot7^3\cdot11\cdot13\cdot17\cdot23\cdot29$.
Its triple cover $3.Fi_{24}'$ embeds into the Monster as class $3A$.
Indeed, on comparing with Table~\ref{t:cent}, we see that the order of $3.Fi_{24}'$ is the size of the centralizer of class $3A$.
Here, the analogue of \eqref{Me8} was again generalized by \cite{HLY} for $Fi_{24}$, the double cover of $Fi_{24}'$.
In particular, the involution classes multiply to only 7 classes which correspond to the affine $\widehat{E_6}$ labels:
\begin{equation}\label{Fie6}
\begin{array}{cccccccccc}
&&&& 1a &&&& \\
&&&& | &&&& \\
&&&& 2a &&&& \\
&&&& | &&&& \\
1a & - & 2b & - & 3a & - & 2a & - & 1a \\
\end{array}
\end{equation}

\subsection{Cusp Numbers}
We now strengthen this correspondence of
{\red
\[
(E_8, \IM), \quad (E_7, \IB), \quad (E_6, Fi_{24}')
\]
}
with a further series of observations.
Recalling our definitions in \S\ref{s:MM}, in \cite{CN}, the cusp number $C$ of the fixing group $\cF(g)$ associated to the class of $g$ is computed.
Moreover, the ``Euler characteristic'' of $\cF(g)$ is also computed; this is the integer $D$ such that $\frac{2 \pi}{3D}$ is the area of the fundamental domain of $\cF(g)$.

For reference, we record the quadruple: (1) class name, (2) cusp number $C$, (3) indicator $D$ for the area of fundamental domain, and (4) normalizer group $\Gamma_0(N)^+ \subset PSL(2,\IR)$ in the notation of Eq.~\eqref{n|h+} for the 194 conjugacy classes of $\IM$, reproduced from Table 2 of \cite{CN}; this is presented in Table \ref{t:cusp}. Note that the $D$ is a multiple of $C$.
Also, we have grouped, as is customary, the Galois conjugates together -- e.g., classes 23A and 23B are combined in 23AB -- because, as aforementioned, they have the same McKay-Thompson series.
For completeness we also tally the occurrences of the cusp numbers within the 172 rational conjugacy classes:
\begin{center}
\begin{tabular}{|c||c|c|c|c|c|c|}\hline
Cusp Number: & 1 & 2 & 3 & 4 & 6 & 8 \\ \hline
Frequency: & 60 & 75 & 12 & 20 & 3 & 2 \\ \hline
\end{tabular}
\end{center}
\comment{
We see that the number of cusp number 1 conjugacy classes is also precisely 60, half of 120.
Moreover, the number of cusp number 2 classes is 75, four times of which, added to 60, again gives 360.
}

We are finally ready to state the first of our key observations, which was in fact made by the second author a number of years ago.
The goal of the remainder of this paper will be to generalize this observation in various contexts.
{\blue
\begin{observation}
For the Monster, we have the following sums for the cusp numbers $C_g$ over the 172 rational conjugacy classes:
\[
\sum_{g} C_g = 360 = 3 \cdot 120\ ,
\qquad
\sum_{g} C_g^2 = 1024 = 2^{10}\ .
\]
\label{360}
\end{observation}
}
We remark that these independent classes have distinct McKay-Thompson series, all except 27A and 27B, which share the same Hauptmodul.

The 360 we recall, from Table~\ref{t:geom}, is thrice 120, which is the number of tritangent planes to Bring's curve. We will generalize this to a wider context of groups shortly.
Furthermore, in analogy to Bring's sextic curve from \eqref{bring}, there is the octavic of Fricke \cite{edge2} of genus 9, the Fermat $[4|1,2,4]$ defined as
\begin{equation}\label{frickeC}
\cF = \{
\sum_i x_i^4 = \sum_i x_i^2 = \sum_i x_i = 0
\} \subset \IP^4 \ .
\end{equation}
The number of tritangent planes on $\cF$ is precisely $2048 = 2 \cdot 1024$, twice the sum of square of the cusps.

Furthermore, in light of \eqref{theta}, the numbers of odd and even theta characteristics on a curve of genus 4, as is the case with Bring's curve, are respectively $2^{4-1}(2^4 - 1) = 120$ and $2^{4-1}(2^4 + 1) = 136$, for a total of $2^8$.
The number of odd $\vartheta$ is precisely the number of bitangents.

\begin{table}[h!!!]
{\tiny
\[
{\hspace{-1cm}
\begin{array}{|llll|}\hline
 \{\text{1A},1,2,1\} & \{\text{2A},1,3,\text{2+}\} &
   \{\text{2B},2,6,\text{2-}\} & \{\text{3A},1,4,\text{3+}\} \\
 \{\text{3B},2,8,\text{3-}\} & \{\text{3C},1,6,\text{3$|$3}\} &
   \{\text{4A},2,6,\text{4+}\} & \{\text{4B},1,6,\text{4$|$2+}\} \\
 \{\text{4C},3,12,\text{4-}\} & \{\text{4D},2,12,\text{4$|$2-}\} &
   \{\text{5A},1,6,\text{5+}\} & \{\text{5B},2,12,\text{5-}\} \\
 \{\text{6A},1,6,\text{6+}\} & \{\text{6B},2,12,\text{6+6}\} &
   \{\text{6C},2,12,\text{6+3}\} & \{\text{6D},2,12,\text{6+2}\} \\
 \{\text{6E},4,24,\text{6-}\} & \{\text{6F},2,18,\text{6$|$3}\} &
   \{\text{7A},1,8,\text{7+}\} & \{\text{7B},2,16,\text{7-}\} \\
 \{\text{8A},2,12,\text{8+}\} & \{\text{8B},2,12,\text{8$|$2+}\} &
   \{\text{8C},1,12,\text{8$|$4}\} & \{\text{8D},4,24,\text{8$|$2-}\} \\
 \{\text{8E},4,24,\text{8-}\} & \{\text{8F},2,24,\text{8$|$4-}\} &
   \{\text{9A},2,12,\text{9+}\} & \{\text{9B},4,24,\text{9-}\} \\
 \{\text{10A},1,9,\text{10+}\} & \{\text{10B},2,18,\text{10+5}\} &
   \{\text{10C},2,18,\text{10+2}\} & \{\text{10D},2,18,\text{10+10}\} \\
 \{\text{10E},4,36,\text{10-}\} & \{\text{11A},1,24,\text{11+}\} &
   \{\text{12A},2,12,\text{12+}\} & \{\text{12B},4,24,\text{12+4}\} \\
 \{\text{12C},1,12,\text{12$|$2+}\} & \{\text{12D},2,18,\text{12$|$3+}\} &
   \{\text{12E},3,24,\text{12+3}\} & \{\text{12F},2,24,\text{12$|$2+6}\}
   \\
 \{\text{12G},2,24,\text{12$|$2+2}\} & \{\text{12H},3,24,\text{12+12}\} &
   \{\text{12I},6,48,\text{12-}\} & \{\text{12J},2,36,\text{12$|$6}\} \\
 \{\text{13A},1,14,\text{13+}\} & \{\text{13B},2,28,\text{13-}\} &
   \{\text{14A},1,12,\text{14+}\} & \{\text{14B},2,24,\text{14+7}\} \\
 \{\text{14C},2,24,\text{14+14}\} & \{\text{15A},1,12,\text{15+}\} &
   \{\text{15B},2,24,\text{15+5}\} & \{\text{15C},2,24,\text{15+15}\} \\
 \{\text{15D},2,36,\text{15$|$3}\} & \{\text{16A},2,24,\text{16$|$2+}\} &
   \{\text{16B},6,48,\text{16-}\} & \{\text{16C},3,24,\text{16+}\} \\
 \{\text{17A},1,18,\text{17+}\} & \{\text{18A},4,36,\text{18+2}\} &
   \{\text{18B},2,18,\text{18+}\} & \{\text{18C},4,36,\text{18+9}\} \\
 \{\text{18D},8,72,\text{18-}\} & \{\text{18E},4,36,\text{18+18}\} &
   \{\text{19A},1,20,\text{19+}\} & \{\text{20A},2,18,\text{20+}\} \\
 \{\text{20B},1,18,\text{20$|$2+}\} & \{\text{20C},4,36,\text{20+4}\} &
   \{\text{20D},2,36,\text{20$|$2+5}\} &
   \{\text{20E},2,36,\text{20$|$2+10}\} \\
 \{\text{20F},3,36,\text{20+20}\} & \{\text{21A},1,16,\text{21+}\} &
   \{\text{21B},2,32,\text{21+3}\} & \{\text{21C},1,24,\text{21$|$3+}\} \\
 \{\text{21D},2,32,\text{21+21}\} & \{\text{22A},1,18,\text{22+}\} &
   \{\text{22B},2,36,\text{22+11}\} & \{\text{23AB},1,24,\text{23+}\} \\
 \{\text{24A},2,24,\text{24$|$2+}\} & \{\text{24B},2,24,\text{24$|$+}\} &
   \{\text{24C},4,48,\text{24$|$+8}\} &
   \{\text{24D},4,48,\text{24$|$2+3}\} \\
 \{\text{24E},2,36,\text{24$|$6+}\} & \{\text{24F},2,48,\text{24$|$4+6}\}
   & \{\text{24G},2,48,\text{24$|$4+2}\} &
   \{\text{24H},4,48,\text{24$|$2+12}\} \\
 \{\text{24I},4,48,\text{24+24}\} & \{\text{24J},2,72,\text{24$|$12}\} &
   \{\text{25A},3,30,\text{25+}\} & \{\text{26A},1,21,\text{26+}\} \\
 \{\text{26B},2,42,\text{26+26}\} & \{\text{27A},3,36,\text{27+}\} &
   \{\text{27B},3,36,\text{27+}\} & \{\text{28A},1,24,\text{28$|$2+}\} \\
 \{\text{28B},2,24,\text{28+}\} & \{\text{28C},3,48,\text{28+7}\} &
   \{\text{28D},2,48,\text{28$|$2+14}\} & \{\text{29A},1,30,\text{29+}\}
   \\
 \{\text{30A},2,36,\text{30+6,10,15}\} & \{\text{30B},1,18,\text{30+}\} &
   \{\text{30C},2,36,\text{30+3,5,15}\} &
   \{\text{30D},2,36,\text{30+5,6,30}\} \\
 \{\text{30E},2,54,\text{30$|$3+10}\} &
   \{\text{30F},2,36,\text{30+2,15,30}\} &
   \{\text{30G},4,72,\text{30+15}\} & \{\text{31AB},1,32,\text{31+}\} \\
 \{\text{32A},4,48,\text{32+}\} & \{\text{32B},4,48,\text{32$|$2+}\} &
   \{\text{33A},2,48,\text{33+11}\} & \{\text{33B},1,24,\text{33+}\} \\
 \{\text{34A},1,27,\text{34+}\} & \{\text{35A},1,24,\text{35+}\} &
   \{\text{35B},2,48,\text{35+35}\} & \{\text{36A},4,36,\text{36+}\} \\
 \{\text{36B},8,32,\text{36+4}\} & \{\text{36C},2,36,\text{36$|$2+}\} &
   \{\text{36D},6,72,\text{36+36}\} & \{\text{38A},1,60,\text{38+}\} \\
 \{\text{39A},1,28,\text{39+}\} & \{\text{39B},1,42,\text{39$|$3+}\} &
   \{\text{39CD},2,56,\text{39+39}\} & \{\text{40A},1,36,\text{40$|$4+}\}
   \\
 \{\text{40B},2,36,\text{40$|$2+}\} & \{\text{40CD},4,72,\text{40$|$2+20}\}
   & \{\text{41A},1,42,\text{41+}\} & \{\text{42A},1,24,\text{42+}\} \\
 \{\text{42B},2,48,\text{42+6,14,21}\} &
   \{\text{42C},2,72,\text{42$|$3+7}\} &
   \{\text{42D},2,48,\text{42$|$3+14,42}\} &
   \{\text{44AB},2,36,\text{44+}\} \\
 \{\text{45A},2,36,\text{45+}\} & \{\text{46AB},2,72,\text{46+23}\} &
   \{\text{46CD},1,36,\text{46+}\} & \{\text{47AB},1,48,\text{47+}\} \\
 \{\text{48A},2,48,\text{48$|$2+}\} & \{\text{50A},3,45,\text{50+}\} &
   \{\text{51A},1,36,\text{51+}\} & \{\text{52A},1,42,\text{52$|$2+}\} \\
 \{\text{52B},2,84,\text{52$|$2+26}\} & \{\text{54A},3,54,\text{54+}\} &
   \{\text{55A},1,36,\text{55+}\} & \{\text{56A},2,48,\text{56+}\} \\
 \{\text{56BC},2,96,\text{56$|$4+14}\} & \{\text{57A},1,60,\text{57$|$3+}\}
   & \{\text{59AB},1,60,\text{59+}\} & \{\text{60A},1,36,\text{60$|$2+}\}
   \\
 \{\text{60B},2,36,\text{60+}\} & \{\text{60C},4,72,\text{60+4,15,60}\} &
   \{\text{60D},3,72,\text{60+12,15,20}\} &
   \{\text{60E},2,72,\text{60$|$2+5,6,30}\} \\
 \{\text{60F},2,108,\text{60$|$6+10}\} & \{\text{62AB},1,48,\text{62+}\} &
   \{\text{66A},1,36,\text{66+}\} & \{\text{66B},2,72,\text{66+6,11,66}\}
   \\
 \{\text{68A},1,54,\text{68$|$2+}\} & \{\text{69AB},1,48,\text{69+}\} &
   \{\text{70A},1,36,\text{70+}\} & \{\text{70B},2,72,\text{70+10,14,35}\}
   \\
 \{\text{71AB},1,72,\text{71+}\} & \{\text{78A},1,42,\text{78+}\} &
   \{\text{78BC},2,84,\text{78+6,26,39}\} &
   \{\text{84A},1,48,\text{84$|$2+}\} \\
 \{\text{84B},2,96,\text{84$|$2+6,14,21}\} &
   \{\text{84C},2,72,\text{84$|$3+}\} & \{\text{87AB},1,60,\text{87+}\} &
   \{\text{88AB},2,72,\text{88$|$2+}\} \\
 \{\text{92AB},2,72,\text{92+}\} & \{\text{93AB},1,96,\text{93$|$3+}\} &
   \{\text{94AB},1,72,\text{94+}\} & \{\text{95AB},1,60,\text{95+}\} \\
 \{\text{104AB},1,84,\text{104$|$4+}\} & \{\text{105A},1,48,\text{105+}\} &
   \{\text{110A},1,54,\text{110+}\} & \{\text{119AB},1,72,\text{119+}\}
\\ \hline
\end{array}
}
\]
}
{\caption
{\sf
The quadruples consisting of (1) class name, (2) cusp number $C$, (3) indicator $D$ for the area of fundamental domain, and (4) normalizer group $\Gamma_0(N)^+ \subset PSL(2,\IR)$ in the notation of Eq.~\eqref{n|h+} for the 172 rational conjugacy classes of Monster group.
}
\label{t:cusp}
}
\end{table}

\subsubsection{Cusp Character}
Let us now consider the full length 194 vector of the cusp numbers, without considering the linear dependencies.
The centralizer $Z(c)$ of each of the 194 conjugacy classes $c$ of the Monster can be found in \cite{atlas} and also in Table 2a of \cite{CN}.
For reference, we give their size (in prime-factorized form) together with the class names, in Table \ref{t:cent}.

\begin{table}[h!!!]
{\scriptsize
\[
{\hspace{-1cm}
\begin{array}{|l|}\hline
 ( {1A} , 2^{46}\cdot 3^{20}\cdot 5^9\cdot 7^6\cdot 11^2\cdot 13^3\cdot
   17\cdot 19\cdot 23\cdot 29\cdot 31\cdot 41\cdot 47\cdot
   59\cdot 71  ) ; \ 
\\
  ( {2A} , 2^{42}\cdot 3^{13}\cdot 5^6\cdot 7^2\cdot 11\cdot 13\cdot
   17\cdot 19\cdot 23\cdot 31\cdot 47  ) ; \
  ( {2B} , 2^{46}\cdot 3^9\cdot 5^4\cdot 7^2\cdot 11\cdot 13\cdot
   23  ) ; \ 
\\
  ( {3A} , 2^{21}\cdot 3^{17}\cdot 5^2\cdot 7^3\cdot 11\cdot 13\cdot
   17\cdot 23\cdot 29  ) ; \ 
  ( {3B} , 2^{14}\cdot 3^{20}\cdot 5^2\cdot 7\cdot 11\cdot 13  ) ; \ 
  ( {3C} , 2^{15}\cdot 3^{11}\cdot 5^3\cdot 7^2\cdot 13\cdot 19\cdot
   31  ) ; \ 
\\
  ( {4A} , 2^{34}\cdot 3^7\cdot 5^3\cdot 7\cdot 11\cdot 23  ) ; \ 
  ( {4B} , 2^{27}\cdot 3^6\cdot 5^2\cdot 7^2\cdot 13\cdot 17  ) ; \ 
  ( {4C} , 2^{34}\cdot 3^4\cdot 5\cdot 7  ) ; \ 
  ( {4D} , 2^{27}\cdot 3^3\cdot 5^2\cdot 7\cdot 13  ) ; \ 
\\
  ( {5A} , 2^{14}\cdot 3^6\cdot 5^7\cdot 7\cdot 11\cdot 19  ) ; \ 
  ( {5B} , 2^8\cdot 3^3\cdot 5^9\cdot 7  ) ; \ 
\\
  ( {6A} , 2^{19}\cdot 3^{10}\cdot 5^2\cdot 7\cdot 11\cdot 13  ) ; \ 
  ( {6B} , 2^{14}\cdot 3^8\cdot 5^2\cdot 7\cdot 11\cdot 13  ) ; \ 
  ( {6C} , 2^{21}\cdot 3^8\cdot 5\cdot 7  ) ; \ 
  ( {6D} , 2^{14}\cdot 3^{13}\cdot 5  ) ; \ 
  ( {6E} , 2^{14}\cdot 3^9\cdot 5  ) ; \ 
  ( {6F} , 2^{15}\cdot 3^5\cdot 5\cdot 7  ) ; \ 
\\
  ( {7A} , 2^{10}\cdot 3^3\cdot 5^2\cdot 7^4\cdot 17  ) ; \ 
  ( {7B} , 2^4\cdot 3^2\cdot 5\cdot 7^6  ) ; \ 
\\
  ( {8A} , 2^{22}\cdot 3^3\cdot 7  ) ; \ 
  ( {8B} , 2^{19}\cdot 3^3\cdot 5\cdot 11  ) ; \ 
  ( {8C} , 2^{14}\cdot 3^3\cdot 5^2\cdot 13  ) ; \ 
  ( {8D} , 2^{19}\cdot 3^2\cdot 5  ) ; \ 
  ( {8E} , 2^{22}\cdot 3  ) ; \ 
  ( {8F} , 2^{14}\cdot 3^3\cdot 7  ) ; \ 
\\
  ( {9A} , 2^6\cdot 3^{11}\cdot 5  ) ; \ 
  ( {9B} , 2^4\cdot 3^{11}  ) ; \ 
  ( {10A} , 2^{11}\cdot 3^2\cdot 5^4\cdot 7\cdot 11  ) ; \ 
  ( {10B} , 2^{14}\cdot 3^2\cdot 5^3  ) ; \ 
  ( {10C} , 2^8\cdot 3\cdot 5^6  ) ; \ 
  ( {10D} , 2^8\cdot 3^3\cdot 5^3\cdot 7  ) ; \ 
  ( {10E} , 2^8\cdot 3\cdot 5^4  ) ; \ 
\\
  ( {11A} , 2^6\cdot 3^3\cdot 5\cdot 11^2  ) ; \ 
  ( {12A} , 2^{15}\cdot 3^6\cdot 5  ) ; \ 
  ( {12B} , 2^{11}\cdot 3^7\cdot 5  ) ; \ 
  ( {12C} , 2^{11}\cdot 3^5\cdot 5\cdot 7  ) ; \ 
  ( {12D} , 2^{11}\cdot 3^4\cdot 7  ) ; \ 
  ( {12E} , 2^{15}\cdot 3^3  ) ; \ 
  ( {12F} , 2^9\cdot 3^3\cdot 5\cdot 7  ) ; \ 
\\
  ( {12G} , 2^9\cdot 3^6  ) ; \ 
  ( {12H} , 2^{11}\cdot 3^3\cdot 5  ) ; \ 
  ( {12I} , 2^{10}\cdot 3^4  ) ; \ 
  ( {12J} , 2^9\cdot 3^2\cdot 5  ) ; \ 
  ( {13A} , 2^4\cdot 3^3\cdot 13^2  ) ; \ 
  ( {13B} , 2^3\cdot 3\cdot 13^3  ) ; \ 
\\
  ( {14A} , 2^9\cdot 3^2\cdot 5\cdot 7^2  ) ; \ 
  ( {14B} , 2^{10}\cdot 3\cdot 7^2  ) ; \ 
  ( {14C} , 2^4\cdot 3^2\cdot 5\cdot 7^2  ) ; \ 
  ( {15A} , 2^6\cdot 3^5\cdot 5^2\cdot 7  ) ; \ 
  ( {15B} , 2^3\cdot 3^6\cdot 5^2  ) ; \ 
  ( {15C} , 2^4\cdot 3^3\cdot 5^2  ) ; \ 
  ( {15D} , 2^3\cdot 3^2\cdot 5^3  ) ; \ 
\\
  ( {16A} , 2^{12}\cdot 3  ) ; \ 
  ( {16B} , 2^{13}  ) ; \ 
  ( {16C} , 2^{13}  ) ; \ 
  ( {17A} , 2^3\cdot 3\cdot 7\cdot 17  ) ; \
  ( {18A} , 2^4\cdot 3^7  ) ; \ 
  ( {18B} , 2^5\cdot 3^6  ) ; \ 
  ( {18C} , 2^6\cdot 3^5  ) ; \ 
\\
  ( {18D} , 2^4\cdot 3^5  ) ; \ 
  ( {18E} , 2^4\cdot 3^5  ) ; \ 
  ( {19A} , 2^2\cdot 3\cdot 5\cdot 19  ) ; \ 
  ( {20A} , 2^{10}\cdot 3\cdot 5^2  ) ; \ 
  ( {20B} , 2^7\cdot 3^2\cdot 5^2  ) ; \ 
  ( {20C} , 2^6\cdot 3\cdot 5^3  ) ; \ 
  ( {20D} , 2^8\cdot 3\cdot 5^2  ) ; \ 
  ( {20E} , 2^4\cdot 3\cdot 5^2  ) ; \ 
\\
  ( {20F} , 2^6\cdot 3\cdot 5  ) ; \ 
  ( {21A} , 2^3\cdot 3^3\cdot 5\cdot 7^2  ) ; \ 
  ( {21B} , 2\cdot 3^2\cdot 7^3  ) ; \ 
  ( {21C} , 2^3\cdot 3^2\cdot 7^2  ) ; \ 
  ( {21D} , 2^3\cdot 3^2\cdot 7  ) ; \ 
  ( {22A} , 2^4\cdot 3\cdot 5\cdot 11  ) ; \ 
  ( {22B} , 2^6\cdot 3\cdot 11  ) ; \ 
\\
  ( {23A} , 2^3\cdot 3\cdot 23  ) ; \ 
  ( {23B} , 2^3\cdot 3\cdot 23  ) ; \ 
  ( {24A} , 2^8\cdot 3^3  ) ; \ 
  ( {24B} , 2^9\cdot 3^2  ) ; \ 
  ( {24C} , 2^7\cdot 3^3  ) ; \ 
  ( {24D} , 2^8\cdot 3^2  ) ; \ 
  ( {24E} , 2^7\cdot 3^2  ) ; \ 
  ( {24F} , 2^5\cdot 3^3  ) ; \ 
  ( {24G} , 2^5\cdot 3^3  ) ; \ 
\\
  ( {24H} , 2^6\cdot 3^2  ) ; \ 
  ( {24I} , 2^7\cdot 3  ) ; \ 
  ( {24J} , 2^5\cdot 3^2  ) ; \ 
  ( {25A} , 2\cdot 5^3  ) ; \ 
  ( {26A} , 2^4\cdot 3\cdot 13  ) ; \ 
  ( {26B} , 2^3\cdot 3\cdot 13  ) ; \ 
  ( {27A} , 2\cdot 3^5  ) ; \ 
  ( {27B} , 3^5  ) ; \ 
\\
  ( {28A} , 2^5\cdot 3\cdot 7^2  ) ; \ 
  ( {28B} , 2^7\cdot 3\cdot 7  ) ; \ 
  ( {28C} , 2^7\cdot 7  ) ; \ 
  ( {28D} , 2^3\cdot 3\cdot 7  ) ; \ 
  ( {29A} , 3\cdot 29  ) ; \ 
\\
  ( {30A} , 2^4\cdot 3^3\cdot 5^2  ) ; \ 
  ( {30B} , 2^5\cdot 3^2\cdot 5^2  ) ; \ 
  ( {30C} , 2^6\cdot 3^2\cdot 5  ) ; \ 
  ( {30D} , 2^3\cdot 3^2\cdot 5^2  ) ; \ 
  ( {30E} , 2^3\cdot 3^2\cdot 5  ) ; \ 
  ( {30F} , 2^4\cdot 3\cdot 5  ) ; \ 
  ( {30G} , 2^4\cdot 3\cdot 5  ) ; \ 
\\
  ( {31A} , 2\cdot 3\cdot 31  ) ; \ 
  ( {31B} , 2\cdot 3\cdot 31  ) ; \ 
  ( {32A} , 2^7  ) ; \ 
  ( {32B} , 2^7  ) ; \ 
  ( {33A} , 2\cdot 3^3\cdot 11  ) ; \ 
  ( {33B} , 2^2\cdot 3^2\cdot 11  ) ; \ 
  ( {34A} , 2^3\cdot 17  ) ; \ 
  ( {35A} , 2^2\cdot 3\cdot 5^2\cdot 7  ) ; \ 
  ( {35B} , 2\cdot 5\cdot 7  ) ; \ 
\\
  ( {36A} , 2^4\cdot 3^4  ) ; \ 
  ( {36B} , 2^3\cdot 3^4  ) ; \ 
  ( {36C} , 2^3\cdot 3^3  ) ; \ 
  ( {36D} , 2^3\cdot 3^2  ) ; \ 
  ( {38A} , 2^2\cdot 19  ) ; \ 
  ( {39A} , 2\cdot 3^3\cdot 13  ) ; \ 
  ( {39B} , 3^2\cdot 13  ) ; \ 
  ( {39C} , 2\cdot 3\cdot 13  ) ; \ 
  ( {39D} , 2\cdot 3\cdot 13  ) ; \
\\ 
  ( {40A} , 2^4\cdot 5^2  ) ; \ 
  ( {40B} , 2^6\cdot 5  ) ; \ 
  ( {40C} , 2^4\cdot 5  ) ; \ 
  ( {40D} , 2^4\cdot 5  ) ; \ 
  ( {41A} , 41  ) ; \ 
  ( {42A} , 2^3\cdot 3^2\cdot 7  ) ; \ 
  ( {42B} , 2^3\cdot 3^2\cdot 7  ) ; \ 
  ( {42C} , 2^3\cdot 3\cdot 7  ) ; \ 
  ( {42D} , 2\cdot 3^2\cdot 7  ) ; \
\\ 
  ( {44A} , 2^5\cdot 11  ) ; \ 
  ( {44B} , 2^5\cdot 11  ) ; \ 
  ( {45A} , 3^3\cdot 5  ) ; \ 
  ( {46A} , 2^3\cdot 23  ) ; \ 
  ( {46B} , 2^3\cdot 23  ) ; \ 
  ( {46C} , 2^2\cdot 23  ) ; \ 
  ( {46D} , 2^2\cdot 23  ) ; \ 
  ( {47A} , 2\cdot 47  ) ; \ 
  ( {47B} , 2\cdot 47  ) ; \ 
\\
  ( {48A} , 2^5\cdot 3  ) ; \ 
  ( {50A} , 2\cdot 5^2  ) ; \ 
  ( {51A} , 3\cdot 17  ) ; \ 
  ( {52A} , 2^3\cdot 13  ) ; \ 
  ( {52B} , 2^2\cdot 13  ) ; \ 
  ( {54A} , 2\cdot 3^3  ) ; \ 
  ( {55A} , 2\cdot 5\cdot 11  ) ; \ 
\\
  ( {56A} , 2^4\cdot 7  ) ; \ 
  ( {56B} , 2^3\cdot 7  ) ; \ 
  ( {56C} , 2^3\cdot 7  ) ; \ 
  ( {57A} , 3\cdot 19  ) ; \ 
  ( {59A} , 59  ) ; \ 
  ( {59B} , 59  ) ; \ 
\\
  ( {60A} , 2^3\cdot 3^2\cdot 5  ) ; \ 
  ( {60B} , 2^4\cdot 3\cdot 5  ) ; \ 
  ( {60C} , 2^3\cdot 3\cdot 5  ) ; \ 
  ( {60D} , 2^3\cdot 3\cdot 5  ) ; \ 
  ( {60E} , 2^2\cdot 3\cdot 5  ) ; \ 
  ( {60F} , 2^2\cdot 3\cdot 5  ) ; \ 
  ( {62A} , 2\cdot 31  ) ; \ 
  ( {62B} , 2\cdot 31  ) ; \ 
\\
  ( {66A} , 2^2\cdot 3\cdot 11  ) ; \ 
  ( {66B} , 2\cdot 3\cdot 11  ) ; \ 
  ( {68A} , 2^2\cdot 17  ) ; \ 
  ( {69A} , 3\cdot 23  ) ; \ 
  ( {69B} , 3\cdot 23  ) ; \ 
  ( {70A} , 2^2\cdot 5\cdot 7  ) ; \ 
  ( {70B} , 2\cdot 5\cdot 7  ) ; \ 
  ( {71A} , 71  ) ; \ 
  ( {71B} , 71  ) ; \ 
\\
  ( {78A} , 2\cdot 3\cdot 13  ) ; \ 
  ( {78B} , 2\cdot 3\cdot 13  ) ; \ 
  ( {78C} , 2\cdot 3\cdot 13  ) ; \ 
  ( {84A} , 2^2\cdot 3\cdot 7  ) ; \ 
  ( {84B} , 2^2\cdot 3\cdot 7  ) ; \ 
  ( {84C} , 2^2\cdot 3\cdot 7  ) ; \ 
  ( {87A} , 3\cdot 29  ) ; \ 
  ( {87B} , 3\cdot 29  ) ; \
\\ 
  ( {88A} , 2^3\cdot 11  ) ; \ 
  ( {88B} , 2^3\cdot 11  ) ; \ 
  ( {92A} , 2^2\cdot 23  ) ; \ 
  ( {92B} , 2^2\cdot 23  ) ; \ 
  ( {93A} , 3\cdot 31  ) ; \ 
  ( {93B} , 3\cdot 31  ) ; \ 
  ( {94A} , 2\cdot 47  ) ; \ 
  ( {94B} , 2\cdot 47  ) ; \ 
  ( {95A} , 5\cdot 19  ) ; \ 
  ( {95B} , 5\cdot 19  ) ; \ 
\\
  ( {104A} , 2^3\cdot 13  ) ; \ 
  ( {104B} , 2^3\cdot 13  ) ; \ 
  ( {105A} , 3\cdot 5\cdot 7  ) ; \ 
  ( {110A} , 2\cdot 5\cdot 11  ) ; \ 
  ( {119A} , 7\cdot 17  ) ; \ 
  ( {119B} , 7\cdot 17  )  \ \\ \hline
\end{array}
}
\]
}
{\caption
{\sf
The size, in prime-factorized form, of the centralizers of each of the 194 conjugacy classes of the Monster, together with their class names.
The centralizer of class $1A$ is the full Monster group.
}\label{t:cent}
}
\end{table}

We see, for example, that the size of the centralizer $Z(1A)$ for the identity class $1A$, is $|M|$. In general, we have
\begin{equation}
|Z(c)| \cdot |c| = |M|
\end{equation}
for each of the 194 conjugacy classes $c$.
Indeed, we have the orthonormality condition for any character table $T_{i\gamma} := \chi_i(c_\gamma)$ of a finite group $G$ with $i=1,2,\ldots,n$ indexing the irreducible representations and $\gamma = 1,2,\ldots,n$ indexing the conjugacy classes. 
The condition states that the weighted table is unitary:
\begin{equation}\label{schur}
U^H \cdot U = \II_{n\times n} \ , \qquad
U_{i \gamma} := T_{i \gamma} |Z(c_\gamma)|^{-\frac12} 
= T_{i \gamma} \sqrt{\frac{|c_\gamma|}{|G|}} \ .
\end{equation}
Less succinctly, the above is customarily presented as the following relations
\begin{align}
\nn
\mbox{Row Orthgonality: } & 
\frac{1}{|G|} \sum\limits_{\gamma = 1}^n \overline{\chi_i(c_\gamma)}
  \chi_j(c_\gamma) |c_\gamma| = \delta_{ij} \ ; 
\\
\mbox{Column Orthgonality: } &
\frac{1}{|G|} \sum\limits_{i = 1}^n \overline{\chi_i(c_\gamma)}
  \chi_i(c_\beta) \sqrt{|c_\gamma||c_\beta|} = \delta_{\gamma \beta} \ . 
\end{align}

Now, consider the list of centralizer sizes $|Z_\gamma| = |G| / |c_\gamma|$ for $\gamma = 1, 2, \ldots, n$. 
This is the character of a reducible representation, which we call the {\it centralizing representation} $R_Z$.
Let $R_Z = \bigoplus_{i=1}^n R_i^{\oplus a_i}$  be expanded into the irreducible representations $R_i$ with coefficients $a_i \in \IZ_{\ge 0}$, so that $\chi(R_Z(c_\gamma)) = \sum\limits_{i=1}^n a_i \chi_i(c_\gamma)$.
We can then use row orthogonality to invert this to obtain
\begin{equation}
a_j = \frac{1}{|G|} \sum\limits_{\gamma=1}^n \frac{|G|}{|c_\gamma|} \chi_j(c_\gamma) |c_\gamma| = \sum\limits_{\gamma=1}^n \chi_j(c_\gamma) \ , \quad
j = 1, 2, \ldots, n \ .
\end{equation}
The sum over the algebraic conjugate representations makes the total sum over the rows of the characters integers, as required.

The above are generalities, which we can certainly apply to the Monster. For instance, the multiplicity coefficients for its centralizing representation begin with $a_j = 194, 203334, 21397838 \ldots$
However, let us now consider the vector of cusp numbers $C_\gamma$ weighted by $|Z(c_\gamma)|$
\begin{equation}
v_\gamma = C_\gamma |Z(c_\gamma)| = C_\gamma \frac{|M|}{|c_\gamma|} \ .
\end{equation}
Is this a character of a representation?

Let us expand as above, i.e., $v_\gamma = \sum\limits_{i=1}^n b_i \chi_i(c_\gamma)$. Inverting using \eqref{schur}, we obtain
\begin{equation}\label{bj}
b_j = \frac{1}{|G|} \sum\limits_{\gamma=1}^n \frac{|G|}{|c_\gamma|} C_\gamma \chi_j(c_\gamma) |c_\gamma| = \sum\limits_{\gamma=1}^n \chi_j(c_\gamma) C_\gamma \ , \quad
j = 1, 2, \ldots, n \ .
\end{equation}
We find the 194 coefficients and see that they are all positive integers! 
Due to their sizes, we present only the first few:
\begin{equation}
b_i = \{
2^3\cdot 7^2, \ 
7^2\cdot 11\cdot 379, \ 
17\cdot 29\cdot 43403, \ 
2063\cdot 409043 \ldots
\}
\end{equation}
That all coefficients are positive integers is non-trivial here because the irrational entries in the character table must conspire to cancel in \eqref{bj}.
It means that the weighted cusps actually correspond to the character of a certain reducible non-trivial representation, which we shall call the {\bf cusp representation}.

\subsection{The Baby and $E_7$ again}
Given that Moonshine has been extended to other groups, even at the very inception of the Monster \cite{CN,queen}, it is only natural to speculate whether other sporadics closely related to $\IM$ give generalizations of Observation \ref{360}, and in particular, ones which touch on the other classical geometries discussed in Table~\ref{t:geom}.
We wish to persist in our $(E_8, \IM)$, $(E_7, \IB)$ and $(E_6, Fi_{24}')$ correspondence. We will see that there indeed is a correspondence, via the cusp numbers, between $E_{8,7,6}$ and the group extensions $\IM, 2.\IB, 3.Fi_{24}'$.

In the context of Moonshine \cite{CN}, we should look at the double cover $2.\IB$, which is associated to class 2A in $\IM$.
Generalized Moonshine for $2.\IB$ has been studied in \cite{baby,Matias}.
Already in the original work of \cite{CN,queen}, the McKay-Thompson series 
was noted to have expansion (with the standard Dedekind eta function $\eta(q)$)
\begin{align}
\nn
T_{2A}(q) & = \left[ \left(\frac{\eta(q)}{\eta(q^2)}\right)^{12} +
2^6 \left(\frac{\eta(q^2)}{\eta(q)}\right)^{12} \right]^2 -104 \\
& =  q^{-1} + 4372q + 96256q^2 + 1240002q^3 + \ldots
\end{align}
which indeed encodes the dimensions of the irreducible representations of $2.\IB$, viz.,\\
$1, 4371, 96255, 1139374 \ldots$

In general, the relevant McKay-Thompson series are cases of the so-called {\bf replicable functions} \cite{FMN,ACMN,mahler}, which have been tabulated comprehensively in \cite{nortonLect,cummins}.
The notation is now standard and is in accord with Norton's database compiled over the years \cite{FMN,nortonLect}.
Of the 616 replicable functions, anything in the form of a number followed by a capital letter is a modular function (of some group between $\Gamma_0(N)$ and $\Gamma_0(N)+$) associated to the matching conjugacy class of $\IM$, i.e., they are monstrous principal moduli. Non-monstrous modular forms are named by a number followed by a small letter, or by a tilde, and then a small letter:
\begin{equation}
\mbox{Monstrous: } nX \ , \qquad
\mbox{other: } nx \quad \mbox{ or } \quad n\sim x \ .
\end{equation}

\begin{table}[h!!!]
{\scriptsize
\[
\hspace{-1cm}
\begin{array}{|rrrrrrrr|} \hline
 \{\text{1a},\text{2A},1\} & \{\text{2a},\text{4$\sim $b},1\} &
   \{\text{2b},\text{2a},1\} & \{\text{2C},\text{4A},2\} &
   \{\text{2d},\text{2B},2\} & \{\text{2e},\text{4C},3\} &
   \{\text{3a},\text{6A},1\} & \{\text{3b},\text{6D},2\} \\
 \{\text{4a},\text{8$\sim $b},1\} & \{\text{4b},\text{4a},1\} &
   \{\text{4c},\text{4B},1\} & \{\text{4d},\text{4C},3\} &
   \{\text{4e},\text{8a},3\} & \{\text{4f},\text{8A},2\} &
   \{\text{4g},\text{8$\sim $d},2\} & \{\text{4h},\text{4D},2\} \\
 \{\text{4i},\text{8B},2\} & \{\text{4j},\text{8E},4\} &
   \{\text{4k},\text{8D},4\} & \{\text{5a},\text{10A},1\} &
   \{\text{5b},\text{10C},2\} & \{\text{6a},\text{12$\sim $d},1\} &
   \{\text{6b},\text{12$\sim $f},2\} & \{\text{6c},\text{6a},1\} \\
 \{\text{6d},\text{6b},1\} & \{\text{6e},\text{12A},2\} &
   \{\text{6f},\text{6C},2\} & \{\text{6g},\text{6c},2\} &
   \{\text{6h},\text{12c},3\} & \{\text{6i},\text{12B},4\} &
   \{\text{6j},\text{6E},4\} & \{\text{6k},\text{12E},3\} \\
 \{\text{6l},\text{12H},3\} & \{\text{6m},\text{12$\sim $h},3\} &
   \{\text{6n},\text{12I},6\} & \{\text{7a},\text{14A},1\} &
   \{\text{8a},\text{16$\sim $a},1\} & \{\text{8b},\text{8a},3\} &
   \{\text{8c},\text{8b},1\} & \{\text{8d},\text{8c},1\} \\
 \{\text{8e},\text{8B},2\} & \{\text{8f},\text{8C},1\} &
   \{\text{8g},\text{8D},4\} & \{\text{8h},\text{16$\sim $d},2\} &
   \{\text{8i},\text{8E},4\} & \{\text{8j},\text{8F},2\} &
   \{\text{8k},\text{16A},2\} & \{\text{8l},\text{16C},3\} \\
 \{\text{8m},\text{16$\sim $e},3\} & \{\text{8n},\text{16a},2\} &
   \{\text{8o},\text{16B},6\} & \{\text{8p},\text{16d},6\} &
   \{\text{9a},\text{18A},4\} & \{\text{9b},\text{18B},2\} &
   \{\text{10a},\text{20$\sim $c},1\} & \{\text{10b},\text{20$\sim $d},2\}
   \\
 \{\text{10c},\text{10a},1\} & \{\text{10d},\text{20A},2\} &
   \{\text{10e},\text{10B},2\} & \{\text{10f},\text{20C},4\} &
   \{\text{10g},\text{10E},4\} & \{\text{10h},\text{20d},3\} &
   \{\text{10i},\text{20F},3\} & \{\text{10j},\text{20$\sim $g},3\} \\
 \{\text{11a},\text{22A},1\} & \{\text{12a},\text{24$\sim $f},1\} &
   \{\text{12b},\text{24$\sim $h},2\} & \{\text{12c},\text{12a},1\} &
   \{\text{12d},\text{12G},2\} & \{\text{12e},\text{12b},1\} &
   \{\text{12f},\text{12C},1\} & \{\text{12g},\text{24a},3\} \\
 \{\text{12h},\text{24$\sim $j},2\} & \{\text{12i},\text{24$\sim $k},2\} &
   \{\text{12j},\text{12E},3\} & \{\text{12k},\text{24b},3\} &
   \{\text{12l},\text{12d},2\} & \{\text{12m},\text{24B},2\} &
   \{\text{12n},\text{24$\sim $m},2\} & \{\text{12o},\text{24c},6\} \\
 \{\text{12p},\text{24A},2\} & \{\text{12q},\text{12I},6\} &
   \{\text{12r},\text{24C},4\} & \{\text{12s},\text{24$\sim $o},4\} &
   \{\text{12t},\text{12F},2\} & \{\text{12u},\text{24h},4\} &
   \{\text{12v},\text{24$\sim $q},4\} & \{\text{12w},\text{24H},4\} \\
 \{\text{12x},\text{24I},4\} & \{\text{12y},\text{24$\sim $r},4\} &
   \{\text{13a},\text{26A},1\} & \{\text{14a},\text{28$\sim $c},1\} &
   \{\text{14b},\text{14a},1\} & \{\text{14c},\text{14c},1\} &
   \{\text{14d},\text{28B},2\} & \{\text{14e},\text{14B},2\} \\
 \{\text{14f},\text{28C},3\} & \{\text{15a},\text{30B},1\} &
   \{\text{15b},\text{30F},2\} & \{\text{16a},\text{16b},3\} &
   \{\text{16b},\text{16c},3\} & \{\text{16c},\text{16B},6\} &
   \{\text{16d},\text{16A},2\} & \{\text{16e},\text{16a},2\} \\
 \{\text{16f},\text{16A},2\} & \{\text{16g},\text{32B},4\} &
   \{\text{16h},\text{32A},4\} & \{\text{16i},\text{32$\sim $e},4\} &
   \{\text{17a},\text{34A},1\} & \{\text{18a},\text{36$\sim $h},4\} &
   \{\text{18b},\text{36$\sim $e},2\} & \{\text{18c},\text{18c},2\} \\
 \{\text{18d},\text{18c},2\} & \{\text{18e},\text{36A},4\} &
   \{\text{18f},\text{18C},4\} & \{\text{18g},\text{36f},6\} &
   \{\text{18h},\text{36D},6\} & \{\text{18i},\text{36$\sim $q},6\} &
   \{\text{19a},\text{38A},1\} & \{\text{20a},\text{40$\sim $c},1\} \\
 \{\text{20b},\text{20a},1\} & \{\text{20c},\text{20c},2\} &
   \{\text{20d},\text{20b},1\} & \{\text{20e},\text{20B},1\} &
   \{\text{20f},\text{40a},3\} & \{\text{20g},\text{40B},2\} &
   \{\text{20h},\text{40$\sim $h},2\} & \{\text{20i},\text{40$\sim $i},2\}
   \\
 \{\text{20j},\text{20E},2\} & \{\text{20k},\text{40C},4\} &
   \{\text{21a},\text{42A},1\} & \{\text{22a},\text{44$\sim $b},1\} &
   \{\text{22b},\text{22a},1\} & \{\text{22c},\text{22a},1\} &
   \{\text{22d},\text{44A},2\} & \{\text{22e},\text{22B},2\} \\
 \{\text{23a},\text{46C},1\} & \{\text{23b},\text{46C},1\} &
   \{\text{24a},\text{48$\sim $c},2\} & \{\text{24b},\text{24d},1\} &
   \{\text{24c},\text{24e},1\} & \{\text{24d},\text{24g},1\} &
   \{\text{24e},\text{24f},1\} & \{\text{24f},\text{24b},3\} \\
 \{\text{24g},\text{24c},6\} & \{\text{24h},\text{24A},2\} &
   \{\text{24i},\text{48$\sim $h},2\} & \{\text{24j},\text{48$\sim $i},2\}
   & \{\text{24k},\text{48A},2\} & \{\text{24l},\text{24H},4\} &
   \{\text{24m},\text{48$\sim $j},3\} & \{\text{24n},\text{48$\sim $k},3\}
   \\
 \{\text{24o},\text{24F},2\} & \{\text{24p},\text{48h},6\} &
   \{\text{25a},\text{50A},3\} & \{\text{26a},\text{52$\sim $c},1\} &
   \{\text{26b},\text{26a},1\} & \{\text{27a},\text{54A},3\} &
   \{\text{28a},\text{56$\sim $d},1\} & \{\text{28b},\text{28A},1\} \\
 \{\text{28c},\text{28C},3\} & \{\text{28d},\text{28a},1\} &
   \{\text{28e},\text{56a},3\} & \{\text{28f},\text{56A},2\} &
   \{\text{28g},\text{56$\sim $g},2\} & \{\text{30a},\text{60$\sim $c},1\}
   & \{\text{30b},\text{60$\sim $l},2\} & \{\text{30c},\text{30a},1\} \\
 \{\text{30d},\text{30d},1\} & \{\text{30e},\text{60B},2\} &
   \{\text{30f},\text{30C},2\} & \{\text{30g},\text{60a},3\} &
   \{\text{30h},\text{60D},3\} & \{\text{30i},\text{60$\sim $m},3\} &
   \{\text{30j},\text{60C},4\} & \{\text{30k},\text{30G},4\} \\
 \{\text{30l},\text{60C},4\} & \{\text{30m},\text{30G},4\} &
   \{\text{31a},\text{62A},1\} & \{\text{31b},\text{62A},1\} &
   \{\text{32a},\text{32B},4\} & \{\text{32b},\text{32B},4\} &
   \{\text{32c},\text{32b},2\} & \{\text{32d},\text{32b},2\} \\
 \{\text{33a},\text{66A},1\} & \{\text{34a},\text{68$\sim $b},1\} &
   \{\text{34b},\text{34a},1\} & \{\text{34c},\text{34a},1\} &
   \{\text{35a},\text{70A},1\} & \{\text{36a},\text{72$\sim $c},2\} &
   \{\text{36b},\text{36C},2\} & \{\text{36c},\text{72a},6\} \\
 \{\text{36d},\text{72$\sim $p},4\} & \{\text{36e},\text{72$\sim $q},4\} &
   \{\text{38a},\text{76$\sim $b},1\} & \{\text{38b},\text{38a},1\} &
   \{\text{38c},\text{38a},1\} & \{\text{39a},\text{78A},1\} &
   \{\text{40a},\text{80$\sim $a},1\} & \{\text{40b},\text{40b},1\} \\
 \{\text{40c},\text{40c},1\} & \{\text{40d},\text{40A},1\} &
   \{\text{40e},\text{80a},2\} & \{\text{40f},\text{80$\sim $e},2\} &
   \{\text{40g},\text{80$\sim $e},2\} & \{\text{42a},\text{84$\sim $e},1\}
   & \{\text{42b},\text{42a},1\} & \{\text{42c},\text{42b},1\} \\
 \{\text{44a},\text{44c},1\} & \{\text{44b},\text{44c},1\} &
   \{\text{46a},\text{92$\sim $b},1\} & \{\text{46b},\text{92$\sim $b},1\}
   & \{\text{46c},\text{92A},2\} & \{\text{46d},\text{46A},2\} &
   \{\text{46e},\text{92A},2\} & \{\text{46f},\text{46A},2\} \\
 \{\text{47a},\text{94A},1\} & \{\text{47b},\text{94A},1\} &
   \{\text{48a},\text{48a},3\} & \{\text{48b},\text{48b},3\} &
   \{\text{50a},\text{100$\sim $c},3\} & \{\text{52a},\text{104$\sim
   $c},1\} & \{\text{54a},\text{108$\sim $g},3\} &
   \{\text{55a},\text{110A},1\} \\
 \{\text{56a},\text{56a},3\} & \{\text{56b},\text{56a},3\} &
   \{\text{60a},\text{120$\sim $d},1\} & \{\text{60b},\text{60b},1\} &
   \{\text{60c},\text{120a},3\} & \{\text{60d},\text{120$\sim $g},2\} &
   \{\text{60e},\text{120$\sim $h},2\} & \{\text{62a},\text{124$\sim
   $b},1\} \\
 \{\text{62b},\text{124$\sim $b},1\} & \{\text{66a},\text{132$\sim $b},1\}
   & \{\text{66b},\text{66a},1\} & \{\text{66c},\text{66a},1\} &
   \{\text{68a},\text{136$\sim $c},1\} & \{\text{70a},\text{140$\sim
   $b},1\} & \{\text{70b},\text{70a},1\} & \{\text{70c},\text{70a},1\} \\
 \{\text{78a},\text{156$\sim $d},1\} & \{\text{84a},\text{168$\sim $c},2\}
   & \{\text{94a},\text{188$\sim $b},1\} & \{\text{94b},\text{188$\sim
   $b},1\} & \{\text{104a},\text{208$\sim $a},1\} &
   \{\text{104b},\text{208$\sim $a},1\} & \{\text{110a},\text{220$\sim
   $b},1\} &  \\
\hline
\end{array}
\]
}
{\caption
{\sf For the 247 classes of the group $2.\IB$, each is a triple $\{mx, nX, c\}$ where $mx$ is the class-name in GAP notation, $nX$ is the identifier for the McKay-Thompson series for the class in the notation of \cite{FMN} and $c$ is the cusp number of the associated modular subgroup.
Of the McKay-Thompson series, 207 are distinct.
}\label{t:baby}
}
\end{table}

Now, there are 247 conjugacy classes of the group $2.\IB$ (for the baby $\IB$ herself, there are 184 conjugacy classes.) and we will name them according to GAP's database, which is also in accord with standard literature \cite{gap}. 
The names are also in the form of - and not to be confused with the McKay-Thompson series - number followed by lower-case letter.
Here, the numbers go from 1 to 110 and the letters go from ``a'' to ``x'' variously.
Using this and the above notation for the McKay-Thompson series for the associated modular subgroup, we can combine the tables of \cite{baby,Matias} and Tables 2 and 3 of \cite{cummins} to obtain all the cusp numbers of these 247 classes.
This is presented in Table~\ref{t:baby} for the reader's convenience.

As with the case of the Monster, we remove the duplicates where different conjugacy classes correspond to the same McKay-Thompson series in Table~\ref{t:baby}, which gives us 207 independent classes over which we can, much as before, sum the cusp numbers as well as their squares.
We arrive at
{\blue
\begin{observation}
For $2.\IB$, we have the following sums for the cusp numbers $C_g$ over the 207 conjugacy classes with distinct McKay-Thompson series:
\[
\sum_{g} C_g(2.\IB) = 448 = 2^3 \cdot 56 \ ,
\qquad
\sum_{g} C_g^2(2.\IB) = 1320 = 2^{3} \cdot 3 \cdot 5 \cdot 11 \ .
\]
\label{448}
\end{observation}
}
Again, examining Table~\ref{t:geom}, the 448 is $2^3$ times 56, the dimension of the fundamental representation of $E_7$.
Likewise, it is a simple (power of 2) factor of 28, which is the number of bitangent lines for $E_7$.

We make two further remarks.
First, there is a total (with repeats) of 106 of classes in $2.\IB$ which are Monstrous (having capital letters) and the sum over cusps $C_g$ is
$
\sum_{g \in \IM} C_g(2.\IB) = 266 = 2 \times 133
$.
We recognize 133 as the complex dimension (and likewise 266 as the real dimension) of $E_7$.
Furthermore, if we only take the rational classes, of which there is a total of 226 (so indeed there are some of these which share the same McKay-Thompson series), the cusp sum becomes simply $2^9 = 512$.

\subsection{Fischer's Group}
Having related the baby to $E_7$ in our context, as discussed in \S\ref{s:e6} and \S\ref{s:e7}, the next natural group to consider is $Fi_{24}'$, the largest of Fischer's sporadic groups.
Now, its triple cover $3.Fi_{24}'$ corresponds to class 3A of the Monster whose McKay-Thompson series is
\begin{align}
\nn
T_{3A}(q) & = \left[ \left(\frac{\eta(q)}{\eta(q^3)}\right)^{6} +
3^3 \left(\frac{\eta(q^3)}{\eta(q)}\right)^{6} \right]^2 -42 \\
& =  q^{-1} + 783q + 8672q^2 + 65367q^3 + \ldots
\end{align}
which indeed encodes the dimensions of the irreducible representations of $3.Fi_{24}'$, viz.,
$1, 8671, 57477 \ldots$

There are 265 conjugacy classes of $3.Fi_{24}'$ in total which again, in standard GAP notation, are labeled as $1a, 2a, \ldots, 105b$.
Amongst these 108 come from $Fi_{24}'$ in an obvious way while the remaining appear as pairs of conjugates under the $\IZ_3$-action; these 108 classes are called {\it essential} in \cite{Matias}.
Of course, some essentials embed into the 265 as singlets and have no Galois orbits of size 3.
Generalized Moonshine for $3.Fi_{24}'$ was studied in \cite{Matias} where all the McKay-Thompson series for the 108 essentials were explicitly constructed.
The orbit classes in the full 265 have McKay-Thompson series being multiplied by $q^{1/3}$ and $q^{2/3}$ and are in some sense not new.

\begin{table}[ht!!!]
{\scriptsize
\[
\hspace{-1.5cm}
\begin{array}{|rrrrr|} \hline
 \{\{\text{1a},\text{3a},\text{3b}\},\text{3A},1\} &
   \{\{\text{2a},\text{6a},\text{6b}\},\text{6A},1\} &
   \{\{\text{2b},\text{6c},\text{6d}\},\text{6C},2\} &
   \{\{\text{3c}\},\text{3C},1\} &
   \{\{\text{3d},\text{3e},\text{3f}\},\text{3B},2\} \\
 \{\{\text{3g},\text{3h},\text{3i}\},\text{9A},2\} &
   \{\{\text{3j}\},\text{9B},4\} & \{\{\text{3k}\},\text{9a},1\} &
   \{\{\text{4a},\text{12a},\text{12b}\},\text{12A},2\} &
   \{\{\text{4b},\text{12c},\text{12d}\},\text{12C},1\} \\
 \{\{\text{4c},\text{12e},\text{12f}\},\text{12E},3\} &
   \{\{\text{5a},\text{15a},\text{15b}\},\text{15A},1\} &
   \{\{\text{6e},\text{6f},\text{6g}\},\text{6d},1\} &
   \{\{\text{6h},\text{6i},\text{6j}\},\text{18$\sim $a},2\} &
   \{\{\text{6k},\text{6l},\text{6m}\},\text{6D},2\} \\
 \{\{\text{6n}\},\text{6F},2\} &
   \{\{\text{6o},\text{6p},\text{6q}\},\text{6E},4\} &
   \{\{\text{6r},\text{6s},\text{6t}\},\text{18B},2\} &
   \{\{\text{6u},\text{6v},\text{6w}\},\text{18A},4\} &
   \{\{\text{6x},\text{6y},\text{6z}\},\text{18E},4\} \\
 \{\{\text{6aa},\text{6ab},\text{6ac}\},\text{18C},4\} &
   \{\{\text{6ad}\},\text{18D},8\} &
   \{\{\text{6ae}\},\text{18e},2\} &
   \{\{\text{7a},\text{21a},\text{21b}\},\text{21A},1\} &
   \{\{\text{7b},\text{21c},\text{21d}\},\text{21B},2\} \\
 \{\{\text{8a},\text{24a},\text{24b}\},\text{24A},2\} &
   \{\{\text{8b},\text{24c},\text{24d}\},\text{24B},2\} &
   \{\{\text{8c},\text{24e},\text{24f}\},\text{24D},4\} &
   \{\{\text{9a}\},\text{9B},4\} &
   \{\{\text{9b},\text{9c},\text{9d}\},\text{9b},1\} \\
 \{\{\text{9e}\},\text{9a},1\} & \{\{\text{9f}\},\text{9B},4\} &
   \{\{\text{9g}\},\text{9c},2\} &
   \{\{\text{9h},\text{9i},\text{9j}\},\text{27A},3\} &
   \{\{\text{10a},\text{30a},\text{30b}\},\text{30B},1\} \\
 \{\{\text{10b},\text{30c},\text{30d}\},\text{30C},2\} &
   \{\{\text{11a},\text{33a},\text{33b}\},\text{33B},1\} &
   \{\{\text{12g}\},\text{12D},2\} &
   \{\{\text{12h},\text{12i},\text{12j}\},\text{12B},4\} &
   \{\{\text{12k}\},\text{12D},2\} \\
 \{\{\text{12l},\text{12m},\text{12n}\},\text{12e},1\} &
   \{\{\text{12o},\text{12p},\text{12q}\},\text{36A},4\} &
   \{\{\text{12r}\},\text{36B},8\} &
   \{\{\text{12s},\text{12t},\text{12u}\},\text{12G},2\} &
   \{\{\text{12v},\text{12w},\text{12x}\},\text{12I},6\} \\
 \{\{\text{12y}\},\text{36b},2\} & \{\{\text{12z}\},\text{36b},2\}
   & \{\{\text{12aa},\text{12ab},\text{12ac}\},\text{36$\sim
   $l},3\} &
   \{\{\text{12ad},\text{12ae},\text{12af}\},\text{36C},2\} &
   \{\{\text{12ag},\text{12ah},\text{12ai}\},\text{36D},6\} \\
 \{\{\text{13a},\text{39a},\text{39b}\},\text{39A},1\} &
   \{\{\text{14a},\text{42a},\text{42b}\},\text{42A},1\} &
   \{\{\text{14b},\text{42c},\text{42d}\},\text{42D},2\} &
   \{\{\text{15c},\text{15d},\text{15e}\},\text{15a},1\} &
   \{\{\text{15f},\text{15g},\text{15h}\},\text{15B},2\} \\
 \{\{\text{15i},\text{15j},\text{15k}\},\text{45A},2\} &
   \{\{\text{16a},\text{48a},\text{48b}\},\text{48A},2\} &
   \{\{\text{17a},\text{51a},\text{51b}\},\text{51A},1\} &
   \{\{\text{18a}\},\text{18e},2\} &
   \{\{\text{18b},\text{18c},\text{18d}\},\text{18h},1\} \\
 \{\{\text{18e}\},\text{18D},8\} &
   \{\{\text{18f},\text{18g},\text{18h}\},\text{18a},1\} &
   \{\{\text{18i}\},\text{18d},2\} &
   \{\{\text{18j},\text{18k},\text{18l}\},\text{54A},3\} &
   \{\{\text{18m},\text{18n},\text{18o}\},\text{54A},3\} \\
 \{\{\text{18p},\text{18q},\text{18r}\},\text{54A},3\} &
   \{\{\text{20a},\text{60a},\text{60b}\},\text{60A},1\} &
   \{\{\text{20b},\text{60c},\text{60d}\},\text{60B},2\} &
   \{\{\text{21e}\},\text{21C},1\} &
   \{\{\text{21f}\},\text{63a},1\} \\
 \{\{\text{21g},\text{21h},\text{21i}\},\text{63$\sim $a},2\} &
   \{\{\text{21j},\text{21k},\text{21l}\},\text{63$\sim $a},2\} &
   \{\{\text{22a},\text{66a},\text{66b}\},\text{66A},1\} &
   \{\{\text{23a},\text{69a},\text{69b}\},\text{69A},1\} &
   \{\{\text{23b},\text{69c},\text{69d}\},\text{69A},1\} \\
 \{\{\text{24g}\},\text{24E},2\} & \{\{\text{24h}\},\text{24E},2\}
   & \{\{\text{24i}\},\text{72b},2\} &
   \{\{\text{24j}\},\text{72b},2\} &
   \{\{\text{24k},\text{24l},\text{24m}\},\text{24C},4\} \\
 \{\{\text{24n},\text{24o},\text{24p}\},\text{72$\sim $r},4\} &
   \{\{\text{24q},\text{24r},\text{24s}\},\text{72$\sim $r},4\} &
   \{\{\text{26a},\text{78a},\text{78b}\},\text{78A},1\} &
   \{\{\text{27a}\},\text{27b},2\} &
   \{\{\text{27b}\},\text{27b},2\} \\
 \{\{\text{27c}\},\text{27b},2\} &
   \{\{\text{28a},\text{84a},\text{84b}\},\text{84A},1\} &
   \{\{\text{29a},\text{87a},\text{87b}\},\text{87A},1\} &
   \{\{\text{29b},\text{87c},\text{87d}\},\text{87A},1\} &
   \{\{\text{30e},\text{30f},\text{30g}\},\text{30b},1\} \\
 \{\{\text{30h},\text{30i},\text{30j}\},\text{90$\sim $a},2\} &
   \{\{\text{33c},\text{33d},\text{33e}\},\text{33A},2\} &
   \{\{\text{33f},\text{33g},\text{33h}\},\text{33A},2\} &
   \{\{\text{35a},\text{105a},\text{105b}\},\text{105A},1\} &
   \{\{\text{36a}\},\text{36b},2\} \\
 \{\{\text{36b}\},\text{36b},2\} & \{\{\text{36c}\},\text{36B},8\}
   & \{\{\text{36d},\text{36e},\text{36f}\},\text{36d},1\} &
   \{\{\text{39c}\},\text{39B},1\} &
   \{\{\text{39d}\},\text{39B},1\} \\
 \{\{\text{39e}\},\text{117a},1\} &
   \{\{\text{39f}\},\text{117a},1\} &
   \{\{\text{42e},\text{42f},\text{42g}\},\text{42c},1\} &
   \{\{\text{42h},\text{42i},\text{42j}\},\text{126$\sim $a},2\} &
   \{\{\text{42k},\text{42l},\text{42m}\},\text{126$\sim $a},2\} \\
 \{\{\text{45a},\text{45b},\text{45c}\},\text{45a},1\} &
   \{\{\text{45d},\text{45e},\text{45f}\},\text{45a},1\} &
   \{\{\text{60e},\text{60f},\text{60g}\},\text{60c},1\} & &
    \\
\hline
\end{array}
\]
}
{\caption
{\sf For the 256 classes of the group $3.Fi_{24}'$, each is a triple $\{\{mx\}, nX, c\}$ where $mx$ is the standard class-name in GAP notation, either as a singlet or as a triplet, $nX$ is the identifier for the McKay-Thompson series for the class in the notation of \cite{FMN} and $c$ is the cusp number of the associated modular subgroup. The triplet $\{mx\}$ is organized according to the $\IZ_3$ Galois orbit of one of the 108 ``essential'' classes of $Fi_{24}'$ which share the same McKay-Thompson series.
}\label{t:Fi24}}
\end{table}

We present, in Table \ref{t:Fi24}, the orbit class structure of the classes, corresponding McKay-Thompson series and the associated cusp number, as a triple in the usual notation
$\{{nx}, mX, c \}$ where $(nx)$ is either a singlet or a triplet of class names depending how an essential class (which is denoted by the first entry) embeds into the full group, $mX$ is Norton's notation for the series and $c$ is the cusp number. Indeed, there will be 108 entries (and on expanding the orbits, the total number of classes is 256).

\comment{
Now, as before, we find the $\IQ$-linearly independent classes (note that this is not necessarily the same as the $\IZ_3$ orbits of the 108 essential classes), we find that there are 165 thereof (including some families of 4 classes):
\begin{equation}
{\scriptsize
\begin{array}{l}
\{\text{1a}\},\{\text{2a}\},\{\text{2b}\},\{\text{3c}\},\{\text{3d}\},\{\text{3g}\},\{\text{3j}\},\{\text{3k}\},\{\text{4a}\},\{\text{4b}\},\{\text{4c}\},\{\text{5a}\},\{\text{6e}\},\{\text{6h}\},\{\text{6k}\},\{\text{6n}\},\{\text{6o}\},\{\text{6r}\},\{\text{6u}\},\{\text{6x}\},\\
\{\text{6aa}\},\{\text{6ad}\},\{\text{6ae}\},\{\text{7a}\},\{\text{7b}\},\{\text{8a}\},\{\text{8b}\},\{\text{8c}\},\{\text{9a}\},\{\text{9b}\},
\{\text{9e}\},\{\text{9f}\},\{\text{9g}\},\{\text{9h}\},\{\text{10a}\},\{\text{10b}\},\{\text{11a}\},\{\text{12g}\},\{\text{12h}\},\{\text{12k}\},\\
\{\text{12l}\},\{\text{12o}\},\{\text{12r}\},\{\text{12s}\},
\{\text{12v}\},\{\text{12y}\},\{\text{12z}\},\{\text{12aa}\},\{\text{12ad}\},\{\text{12ag}\},\{\text{13a}\},\{\text{14a}\},\{\text{14b}\},\{\text{15c}\},\{\text{15f}\},\{\text{15i}\},
\{\text{16a}\},\{\text{17a}\},\\
\{\text{18a}\},\{\text{18b}\},\{\text{18e}\},\{\text{18f}\},\{\text{18i}\},\{\text{18j}\},\{\text{20a}\},\{\text{20b}\},\{\text{21e}\},\{\text{21f}\},\{\text{22a}\},
\{\text{24g}\},\{\text{24h}\},\{\text{24i}\},\{\text{24j}\},\{\text{24k}\},\{\text{26a}\},\\
\{\text{27a}\},\{\text{27b}\},\{\text{27c}\},\{\text{28a}\},\{\text{30e}\},\{\text{30h}\},\{\text{35a}\},
\{\text{36a}\},\{\text{36b}\},\{\text{36c}\},\{\text{36d}\},\{\text{42e}\},\{\text{60e}\},\{\text{3a},\text{3b}\},\{\text{6a},\text{6b}\},\{\text{6c},\text{6d}\},\{\text{3e},\text{3f}\},\{\text{3h},\text{3i}\},\\
\{\text{12a},\text{12b}\},\{\text{12c},\text{12d}\},\{\text{12e},\text{12f}\},\{\text{15a},\text{15b}\},\{\text{6f},\text{6g}\},\{\text{6i},\text{6j}\},\{\text{6l},\text{6m}\},\{\text{6p},\text{6q}\},\{\text{6s},\text{6t}\},
\{\text{6v},\text{6w}\},\{\text{6y},\text{6z}\},\\
\{\text{6ab},\text{6ac}\},\{\text{21a},\text{21b}\},\{\text{21c},\text{21d}\},\{\text{24a},\text{24b}\},\{\text{24c},\text{24d}\},\{\text{24e},\text{24f}\},\{\text{9c},\text{9d}\},\{\text{9i},\text{9j}\},\{\text{30a},\text{30b}\},\{\text{30c},\text{30d}\},\{\text{33a},\text{33b}\},\\
\{\text{12i},\text{12j}\},\{\text{12m},\text{12n}\},\{\text{12p},\text{12q}\},\{\text{12t},\text{12u}\},\{\text{12w},\text{12x}\},\{\text{12ab},\text{12ac}\},\{\text{12ae},\text{12af}\},\{\text{12ah},\text{12ai}\},\{\text{39a},\text{39b}\},\{\text{42a},\text{42b}\},\{\text{42c},\text{42d}\},\\
\{\text{15d},\text{15e}\},\{\text{15g},\text{15h}\},\{\text{15j},\text{15k}\},\{\text{48a},\text{48b}\},\{\text{51a},\text{51b}\},\{\text{18c},\text{18d}\},\{\text{18g},\text{18h}\},\{\text{18k},\text{18l}\},\{\text{18m},\text{18p}\},\{\text{18n},\text{18r}\},\{\text{18o},\text{18q}\},\\
\{\text{60a},\text{60b}\},\{\text{60c},\text{60d}\},\{\text{21g},\text{21j}\},\{\text{66a},\text{66b}\},\{\text{23a},\text{23b}\},\{\text{24l},\text{24m}\},\{\text{24n},\text{24q}\},\{\text{78a},\text{78b}\},\{\text{84a},\text{84b}\},\{\text{29a},\text{29b}\},\{\text{30f},\text{30g}\},\\
\{\text{30i},\text{30j}\},\{\text{33c},\text{33f}\},\{\text{105a},\text{105b}\},\{\text{36e},\text{36f}\},\{\text{39c},\text{39d}\},\{\text{39e},\text{39f}\},\{\text{42f},\text{42g}\},\{\text{42h},\text{42k}\},\{\text{45a},\text{45d}\},\{\text{60f},\text{60g}\},\\
\{\text{21h},\text{21i},\text{21k},\text{21l}\},\{\text{69a},\text{69b},\text{69c},\text{69d}\},\{\text{24o},\text{24p},\text{24r},\text{24s}\},\{\text{87a},\text{87b},\text{87c},\text{87d}\},\{\text{33d},\text{33e},\text{33g},\text{33h}\},\{\text{42i},\text{42j},\text{42l},\text{42m}\},\{\text{45b},\text{45c},\text{45e},\text{45f}\}
\end{array}
}
\end{equation}

From these we readily make the following
{\blue
\begin{observation}
For $3.Fi_{24}'$, we have the following sums for the cusp numbers $C_g$ over the 165 rational conjugacy classes:
\[
\sum_{g} C_g(3.Fi_{24}') = 360 = 3 \cdot 120 \ ,
\qquad
\sum_{g} C_g^2(3.Fi_{24}') = 1146= 2 \cdot 3 \cdot 191 \ .
\]
If we consider the size-4 Galois families (which do not occur for the Monster) as 4 independent single classes, there are 186  $\IQ$-independent conjugacy classes, and the sums become
\begin{equation}
\sum_{g} C_g(3.Fi_{24}') = 399 = 19 \cdot 27 \ ,
\qquad
\sum_{g} C_g^2(3.Fi_{24}') = 1239 = 3 \cdot 7 \cdot 59 \ .
\end{equation}
If we summed over only over the 108 ``essential classes'' which come from $fi_24}'$, then we have that
\begin{equation}
\sum_{g} C_g(3.Fi_{24}') = 240 = 2 \cdot 120 \ ,
\qquad
\sum_{g} C_g^2(3.Fi_{24}') = 802= 2 \cdot 401 \ .
\end{equation}
\end{observation}
}

Now, as before, we extract the classes with unique McKay-Thompson series.
Of the 108 essentials, we see that 83 are distinct, therefrom, likewise in the full group, there will be 213 out of the 256.
However, we are confronted, for the first time, with {\it irrational McKay-Thompson series} due to the multiplication of $q^{\pm 1/3}$, which we could either interpret as being new or not, and we will make the sum in both cases for comparison:
{\blue
\begin{observation}
For $3.Fi_{24}'$, we have the following sums for the cusp numbers $C_g$ over the 213 rational conjugacy classes which have distinct McKay-Thompson series:
\[
\sum_{g} C_g(3.Fi_{24}') = 440 = 2^3 \cdot 5 \cdot 11 =
2^3 \cdot (2 \cdot 27 + 1) \ ,
\qquad
\sum_{g} C_g^2(3.Fi_{24}') =  1290 = 2\cdot3\cdot5\cdot43 \ .
\]
\label{440}
\end{observation}
}
Had we not considered the irrational McKay-Thompson series as distinct but the same as the essential ones, we are then effectively working over the group $Fi_{24}'$, in which case, we have the cusp sums over the 83 distinct classes being
$\sum_{g} C_g(Fi_{24}') = 176 = 2^4 \cdot 11$ and
$\sum_{g} C_g^2(3.Fi_{24}') = 554 = 2 \cdot 277$.

\subsection{Conway's Group}
Other than the Baby, perhaps the closest sporadic groups to the Monster is Conway's group $Co_1$, which is associated to class 2B of $\IM$ by a cover of order $2^{1+24}$.
The sporadic simple group $Co_1$, of order 
$2^{21}\cdot 3^9\cdot 5^4\cdot 7^2\cdot 11\cdot 13\cdot 23$ is itself the $\IZ_2$-quotient of the non-simple group $Co_0$, which is the automorphism group of the famous Leech lattice \footnote{
In analogy to \eqref{thetaE8}, the theta-series for the Leech lattice $\Lambda$ is $\theta_\Lambda(q) = \sum\limits_{x \in \Lambda} q^{|x|^2/2} =\sum\limits_{m=0}^\infty \frac{65520}{691} \left(\sigma_{11}(m) - \tau(m) \right) q^{2m}
= \frac{65520}{691} \left(\sum\limits_{m=0}^\infty \frac{m^{11}q^{2m}}{1 - q^{2m}} - \Delta(q^2)\right)$.
}
, the unique even self-dual lattice in 24-dimensions.
It is worth recalling that $Co_0 = 2.Co_1$ has 167 conjugacy classes while $Co_1$ has 101.

After \cite{queen}, there has been a host of activity to study \cite{koike,kondo,scheit,DC} Moonshine for $Co_0$ as well as $Co_1$, of which we will employ the most recent results in the last reference.
As far back as the earliest results of \cite{koike}, it was realized that the McKay-Thompson series are given explicitly as products and quotients of Dedekind eta-functions whose arguments are appropriate powers of the nome $q$; these are so-called eta-quotients \cite{DKM,martin,He:2013lha}.
For the Mathieu group $M_{24}$ over which there has been extensive activity \cite{Cheng:2013wca,DGO}, all the McKay-Thompson series can be written entirely as eta-products.

The invariance groups 
\footnote{
Along a parallel vein, we can examine the second part of \cite{scheit}, where the relevant results for the square-free case for $Co_0$ are presented in the table under the section entitled ``Genus 0 groups'' at the end.
Importantly, the moonshine (modular) group where $N$ is square-free and in the notation of \cite{CN} is given in column 2, for each of the classes of $Co_0$ in comparison with those of the Monster.
There is a total of 41  square-free conjugacy classes of $Co_0$, corresponding to the following classes of the Monster:
\[
\begin{array}{c}
2B,\ 3B,\ 5B,\ 6B,\ 6C,\ 6D,\ 6E,\ 6E,\ 6E,\ 7B,\  10D,\ 10B,\  10C,\ 10E,\  10E,\ 10E,\ 13B,\  \\
14C,\ 14B,\ 15C,\ 15B,\ 21D,\ 21B,\ 22B,\ 26B,\ 30F,\ 30D,\ 30A,\ 30C,\ 30G,\ 30G,\ 30G,\ \\
33A,\ 35B,\  39CD,\  42D,\  42B,\  46AB,\  66B,\  70B,\  78BC \ .
\end{array}
\]
In the above, other than direct reference to the Table, there are 6 which have no moonshine group labeled explicitly, corresponding to the eta-quotients
\[
\frac{2^83^4}{1^46^8}, \
\frac{2^33^9}{1^36^9}, \
\frac{2^45^2}{1^210^4}, \
\frac{2^15^5}{1^110^5}, \
\frac{1^16^210^215}{2^23^15^130^2}, \
\frac{3.5}{2.30} \ ,
\]
where the notation 
$n^a$ means $\eta(q^{n})^{a}$.
Since we can readily find the $q$-expansions for these, we determine which McKay-Thompson series they are and thus the corresponding class and moonshine group from Table \ref{t:cusp}.
We find that these are, respectively, the series for the classes $6E,6E,10E,10E,30G,30G$.
Checking against the cusp numbers for the Monstrous classes in the above, we readily find that the sum over the cusp numbers is 100 and that of their squares, 272.
Removing repeats, such as precisely the 6 classes of the eta-quotients above , there are 35 classes and the cusp sum now becomes 76 and the square sum, 176.
}
of the classes of $Co_0$ and $Co_1$ are tabulated in the Appendix of \cite{DC}, in the original notation of \cite{CN}.
Unsurprisingly, there are many cases which are {\it not} Monstrous McKay-Thompson series, viz., those which are of the form $nX$ with capital ``X'' in Norton's notation.
Unfortunately, the new ones are not given in the standard $nx$ and $n\sim x$ notation as in our previous cases so extracting their cusp numbers is not immediate.
We leave this exercise to the full study of all the cusp-sums for all the sporadic groups with Moonshine to a future work.
For now, we remark that for those classes of the Conway group which {\it do} have Monstrous McKay-Thompson series, the sums are $\sum_{g} C_g = 165 = 3\cdot5\cdot 11$ and $\sum_{g} C_g^2 = 615 = 3\cdot5\cdot 41$.
In fact, \cite{DC} does more, and lists certain twisted McKay-Thompson series for $Co_0$ and $Co_1$, which are, in fact, all Monstrous (cf~Table 2 in cit.~ibid.), in which case we have 80 distinct (Monstrous) McKay-Thompson series and checking their associated cusp numbers \footnote{
We are grateful to John Duncan for telling us that he has recently computed the full cusp sum for $Co_0$ and the number is 480.
} we have that
\begin{equation}
\sum_{g} C_g(Co_0) = 224 = 2^5 \cdot 7 \ ,
\qquad
\sum_{g} C_g^2(Co_0) = 770 = 2 \cdot5\cdot7\cdot11\ .
\end{equation}

\subsection{Genus Zero}
Finally, it is worth returning to the genus zero principle congruence groups themselves, independently of any particular groups.
Consider the 15 (not necessarily prime) numbers $N$ in \eqref{g=0N} for which the genus of $\Gamma(N)$ is zero.
We can first use the formula \eqref{cuspGamma0N} to obtain their cusp numbers, giving us
\begin{equation}\label{cg0}
\begin{array}{c||c|c|c|c|c|c|c|c|c|c|c|c|c|c|c}
N & 1 & 2 & 3 & 4 & 5 & 6 & 7 & 8 & 9 & 10 & 12 & 13 & 16 & 18 & 25 \\
\hline
c(N) & 1 & 2 & 2 & 3 & 2 & 4 & 2 & 4 & 4 & 4 & 6 & 2 & 6 & 8 & 6 \\
\end{array}
\end{equation}
{\blue
\begin{observation}
For the 15 genus 0 groups $\Gamma_0(N)$ above,
the sum over their cusp numbers is 56, and the sum over their squares is 266.
\end{observation}
}
These are respectively the dimension of the fundamental representation and the real dimension of $E_7$.

We can also look up the cusp numbers for the full normalizer $\Gamma_0(N)^+$ from the entries for $N+$ in the table in \cite{CN} (in the same reference, \eqref{cg0} would be denoted as simply $N-$); they are
$\{
1, 1, 1, 2, 1, 1, 1, 2, 2, 1, 2, 1, 3, 2, 3
\}$ . The sum gives 24, and the sum of squares, 46.

\section{The Horrocks-Mumford Bundle: A Digression}\label{s:HM}\setall
Having addressed Conway's group in our cusp-sporadic correspondence, let us conclude with a parting digression on another context in which the Conway group and the curves of Bring and Fricke arise in relation to a classical geometric object.

Let us investigate the problem of vector bundles on projective spaces, a central subject in algebraic geometry.
Horrocks and Mumford famously constructed their rank 2 indecomposable bundle on $\IP^4$, which is the only known such an example \cite{HM}.
An excellent account, with historical context, is given in \cite{Hulek}.
From this later reference we summarize the following key points about vector bundles on $\IP^n$ (by vector bundles we henceforth mean algebraic, holomorphic, complex vector bundles over projective varieties)
\begin{itemize}
\item[\underline{$n=1$}] 
Grothendieck's theorem \cite{groP1} guarantees that any vector bundle $E$ of rank $r$ on $\IP^1$ splits completely into a direct sum of line bundles as $E = \bigoplus_{i=1}^r \cO_{\IP^1}(a_i)$.

\item[\underline{$n=2$}] 
Wu's theorem \cite{wuP2} states that isomorphism classes of $\IC^2$-vector bundles $E$ on $\IP^2$ are classified by the first and second Chern classes.

\item[\underline{$n=3$}] 
Atiyah-Rees \cite{ARP3} proves that every $\IC^2$-bundle on $\IP^3$ admits an algebraic structure.

\item[\underline{$n\ge4$}] 
Grauert-Schneider \cite{GSP4} conjecture that every unstable rank 2 bundle on $\IP^{n \ge 4}$ splits.

\item[\underline{$n\ge6$}]
Hartshorne \cite{HaP6} conjectures that every rank 2 bundle on $\IP^{n \ge 6}$ splits (this in particular implies that smooth projective varieties $X \subset \IP^n$ for $dim(X) > \frac23 n$ is complete intersection. 
\end{itemize}

One can see that the less is known the higher the dimension of the projective space. In dimension 4, the bundle of Horrocks-Mumford \cite{HM} is essentially the only non-trivial rank 2 bundle \cite{HMunique}, which fits well within our realm of exceptional/sporadic objects.
There are several equivalent constructions \cite{Hulek} and we will follow the so-called {\it monad} construction (cf.~a very explicit computational-geometric description in \cite{HMm2}).
Consider the (non-exact) complex of vector bundles
\begin{equation}\label{HMseq}
0 \longrightarrow
\cO_{\IP^4}(2)^{\oplus 5}
\stackrel{p}{\longrightarrow}
(\bigwedge^2 T )^{\oplus 2}
\stackrel{q}{\longrightarrow}
\cO_{\IP^4}(3)^{\oplus 5}
\longrightarrow
0
\end{equation}
where $T$ is the tangent bundle on $\IP^4$ and $p$ and $q$ are respectively injective and surjective maps of bundles which will be specified in \S\ref{s:HMdetail}.
The Horrocks-Mumford bundle is simply the cohomology of the above complex:
\begin{equation}\label{HMdef}
F_{HM} := \ker(q) / \im(p) \ .
\end{equation}
The total Chern class of the bundle, in terms of the hyperplane class $h$ of $\IP^4$, is
\begin{equation}
c(F_{HM}) = 1 + 5 H + 10 H^2 \ .
\end{equation}
Subsequently, one can use Riemann-Roch to obtain the Hilbert polynomial as
\begin{equation}
\chi(F_{HM} \otimes \cO_{\IP^4}(n H)) = 
\frac{1}{12} n^4 + \frac53 n^3 + \frac{125}{12}n^2 + \frac{125}{6}n + 2 \ .
\end{equation}
In fact, from \S4 of \cite{HM}, we can obtain the Hilbert Series as the generating function of global sections as
\begin{equation}
\sum\limits_{n=0}^\infty h^0(F_{HM}(nH)) t^n =
4 + \frac{1}{12}\sum\limits_{n=6}^\infty (n^2-1)(n^2-24) t^{n-5} =
\frac{4+15t-35t^2+20t^3-2t^5}{(1-t)^5} \ .
\end{equation}
We point out that while in the original paper \cite{HM}, $F_{HM}$ was defined with \eqref{HMseq}, often in later literature \cite{HMm2,HMunique} the dual twisted by $\cO(2)$ is defined as the Horrocks-Mumford bundle, i.e., as the cohomology of the complex
\begin{equation}
0 \longrightarrow
\cO_{\IP^4}(-1)^{\oplus 5}
\stackrel{p'}{\longrightarrow}
(\bigwedge^2 T^* \otimes \cO_{\IP^4}(2))^{\oplus 2}
\stackrel{q'}{\longrightarrow}
\cO_{\IP^4}^{\oplus 5}
\longrightarrow
0
\end{equation}
so that $F' = \ker(q') / \im(p')$, where $T^*$ is now the cotangent bundle of $\IP^4$.
In this definition, $c(F') = 1 - H + 4H^2$.

\subsection{Symmetry Groups}
It was shown in \cite{HM} that the group $G_{HM}$ of symmetries \footnote{
We thank Igor Dolgachev for suggesting to look directly at the polarized Abelian surface in $\IP^4$ of degree 10, on which $G_{HM}$ acts naturally.
} on this bundle is the normalizer of the Heisenberg group $\cH(5)$ within $SL(5; \IQ(\omega_5))$.
Here, we recall that $SL(5; \IQ(\omega_5))$ is the special linear group in dimension 5 defined over the cyclotomic field of $\IQ$ extended by the primitive 5-th root $\omega_5$ of unity.
Furthermore, $\cH(5)$ is the order 125 extra-special group which is the non-Abelian central extension of the Abelian group $\IZ_5 \times \IZ_5$:
\begin{equation}
0 \to \IZ_5 \to \cH(5) \to \IZ_5 \times \IZ_5 \to 0 \ .
\end{equation}
In fact, $G_{HM}$ is a semi-direct product of $\cH(5)$ with the binary icosahedral group $\widehat{E_8} = SL(2,5)$, of order 120.
The Heisenberg group itself acts on the standard basis $e_i$ of $\IC^5$ as:
\begin{equation}
\cH(5) = \left<\sigma, \tau\right> \acts \IC^ \ : \qquad
\sigma :  e_i \longrightarrow e_{i+1} \ , \quad
\tau : e_i \longrightarrow \omega_5^i e_i \ ,
\end{equation}
with subscripts on coordinates defined modulo 5.

The generators \footnote{
We can also represent $G_{HM}$ as a permutation group, acting on 150 elements, in cycle notation, the generators $f_1$ and $f_2$ are
\[
{\tiny
\begin{array}{l}
f_1 = \\
(1,42,15,16,131)(2,43,11,17,132)(3,44,12,18,133)(4,45,13,19,134)\\
(5,41,14,20,135)(6,37,22,88,77)(7,38,23,89,78)(8,39,24,90,79)(9,40,25,86,80)\\
(10,36,21,87,76)(26,141,122,107,67)(27,142,123,108,68)(28,143,124,109,69)(29,
    144,125,110,70)\\
(30,145,121,106,66)(31,92,50,99,129)(32,93,46,100,130)(33,94,47,96,126)(34,95,48,97,127)\\
(35,91,49,98,128)(51,146,71,119,60)(52,147,72,120,56)(53,148,73,116,57)(54,149,74,117,58)\\
(55,150,75,118,59)(61,81,101,138,113)(62,82,102,139,114)(63,83,103,140,115)
(64,84,104,136,111)(65,85,105,137,112) \ ,  \\
f_2 = \\
(1,9,73,36,20)(2,10,74,37,16)(3,6,75,38,17)(4,7,71,39,18)(5,8,72,40,19)\\
(11,82,94,113,53)(12,83,95,114,54)(13,84,91,115,55)(14,85,92,111,51)(15,81,93,112,52)\\
(21,141,123,148,96)(22,142,124,149,97)(23,143,125,150,98)(24,144,121,146,99)(25,145,122,147,100)\\
(26,28,30,27, 29)(31,35,34,33,32)(41,58,64,129,140)(42,59,65,130,136)(43,60,61,126, 137)\\
(44,56,62,127,138)(45,57,63,128,139)(46,118,106,67,77)(47,119,107,68,78)(48,120,108,69,79)\\
(49,116,109,70,80)(50,117,110,66,76)(86,87,88,89,90)(101,103,105,102,104)(131,132,133,134,135)\\
\end{array}
}
\]
}
and presentation of $G_{HM}$ can be computed using \cite{gap} (we record these here because most literature are not explicit about these).
The natural 5-dimensional complex faithful representation gives $G_{HM}$ as a 2-generated group:
\begin{align}
\nn
G_{HM}& =
\left<
{\scriptsize
\left(
\begin{array}{ccccc}
 \omega_5^3 & 0 & 0 & 0 & 0 \\
 0 & \omega_5^2 & 0 & 0 & 0 \\
 0 & 0 & \omega_5^4 & 0 & 0 \\
 0 & 0 & 0 & \omega_5^4 & 0 \\
 0 & 0 & 0 & 0 & \omega_5^2 \\
\end{array}
\right)
} \ ,
\right.
\\
&
\left.
{\tiny
\frac15
\left(
\begin{array}{ccccc}
 -2 \omega_5^3-\omega_5^2-2 \omega_5 & -2 \omega_5^4-\omega_5^3-2 \omega_5^2 & -2 \omega_5^4-\omega_5^3-2 \omega_5^2 & -2 \omega_5^3-\omega_5^2-2
   \omega_5 & -\omega_5^4+\omega_5^3+\omega_5^2-\omega_5 \\
 -2 \omega_5^4-\omega_5^3-2 \omega_5^2 & -\omega_5^4+\omega_5^3+\omega_5^2-\omega_5 & 2 \omega_5^4+2 \omega_5^3+\omega_5 & 2 \omega_5^4+2
   \omega_5^3+\omega_5 & -\omega_5^4+\omega_5^3+\omega_5^2-\omega_5 \\
 -2 \omega_5^4-\omega_5^3-2 \omega_5^2 & 2 \omega_5^4+2 \omega_5^3+\omega_5 & -2 \omega_5^4-\omega_5^3-2 \omega_5^2 & \omega_5^4+2 \omega_5^2+2 \omega_5
   & \omega_5^4+2 \omega_5^2+2 \omega_5 \\
 -2 \omega_5^3-\omega_5^2-2 \omega_5 & 2 \omega_5^4+2 \omega_5^3+\omega_5 & \omega_5^4+2 \omega_5^2+2 \omega_5 & 2 \omega_5^4+2 \omega_5^3+\omega_5 & -2
   \omega_5^3-\omega_5^2-2 \omega_5 \\
 -\omega_5^4+\omega_5^3+\omega_5^2-\omega_5 & -\omega_5^4+\omega_5^3+\omega_5^2-\omega_5 & \omega_5^4+2 \omega_5^2+2 \omega_5 & -2 \omega_5^3-\omega_5^2-2
   \omega_5 & \omega_5^4+2 \omega_5^2+2 \omega_5 \\
\end{array}
\right)
}
\right>
\end{align}

We see that while both semi-direct factors require 3 generators, the full group needs only 2.
Calling the two generators above as $f_1$ and $f_2$, the presentation of the group is simply
\begin{align}
\nn
G_{HM} &:= \left<
f_1^5, \ f_2^5, \ f_2f_1f_2f_1f_2f_1^{-1}f_2^{-1}f_1^{-1}f_2^{-1}f_1^{-1}f_2^{-1}f_1, \ \right.
\\
&
\left.
f_2f_1f_2f_1^2f_2^{-1}f_1f_2^{-1}f_1f_2^{-1}f_1^2, \ 
f_2f_1^{-2}f_2f_1^2f_2f_1^2f_2f_1^{-2}f_2f_1^{-1}
\right> \ .
\end{align}
For reference, we include the character table -- both the ordinary and the rational, as well as the modular versions -- of $G_{HM}$ in Appendix \ref{ap:ct}.
It is interesting to note that the {\it modular} Brauer character table, defined over $\IF_5$ has 5 irreducible representations and 26 conjugacy classes, most of which has 0 Brauer character, except 5. This gives essentially a $5 \times 5$ character table
\begin{equation}
\mbox{Character}_{\IF_5}(G_{HM}) = 
{\scriptsize
\left[
\begin{array}{ccccc}
 1 & 1 & 1 & 1 & 1\\
 2 & -2 & -1 & 0 & 1\\
 3 & 3 & 0 & -1 & 0\\
 4 & -4 & 1 & 0 & -1\\
 5 & 5 & -1 & 1 & -1\\
\end{array}
\right] } \ ,
\end{equation}
and thus it behooves us to consider the McKay quiver \cite{mckay}.
Now, because $G_{HM}$ naturally embeds into $SL(2,\IZ)$, it is expedient to take the fundamental {\bf 2} representation, as in the ADE case of \eqref{ADE}, and decompose ${\bf 2} \otimes {\bf r}_i = \bigoplus_{j=1}^5 a_{ij} {\bf r}_j$.
We readily find by checking, for example, decompositions such as ${\bf 2} \otimes {\bf 2} = {\bf 1} \oplus {\bf 3}$ whose character is $(4,4,1,0,1) = (1,1,1,1,1) + (3,3,0,-1,0)$, that the adjacency matrix $a_{ij}$ and the accompanying McKay quiver as follows:
\begin{equation}
a_{ij} = {\scriptsize 
\left(
\begin{array}{ccccc}
 0 & 1 & 0 & 0 & 0 \\
 1 & 0 & 1 & 0 & 0 \\
 0 & 1 & 0 & 1 & 0 \\
 0 & 0 & 1 & 0 & 1 \\
 0 & 1 & 0 & 2 & 0 \\
\end{array}
\right)}
\qquad
\begin{array}{l}\includegraphics[trim=0mm 0mm 0mm 0mm, width=3.0in]{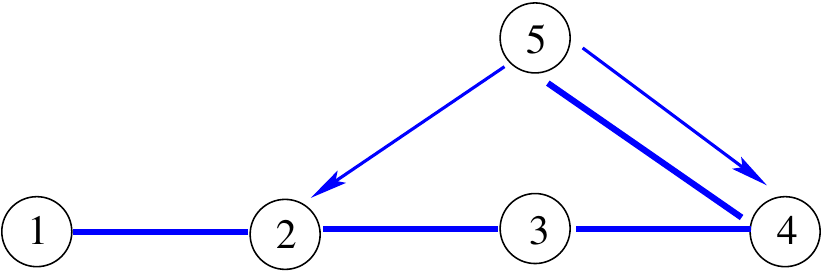}\end{array}
\end{equation}
In the quiver, every undirected edge is a pair of arrows in opposing directions.
The eigenvalues of this adjacency matrix are $\pm2, \pm1, 0$.


\subsection{The HM Quintic Calabi-Yau Threefold}
The space of $\cH(5)$-invariant quintics arose from \cite{HM} as invariant sections of the bundle $\cO_{\IP^4}(5)$ and has been widely studied since, especially in the context of heterotic string compactifications \cite{Braun:2009mb,Anderson:2009mh}.
We recall that this is a 6-dimensional space of quintics in the projective coordinates $[x_0:\ldots:x_4]$ of $\IP^4$,
\begin{align}
\nn
H^0(\IP^4, \cO_{\IP^4}(5))^{\cH(5)} = & \mbox{span}\left<
\sum_i x_i^5 \ , \ \
\sum_i x_i^3 x_{i+1} x_{i+4} \ , \ \
\sum_i x_i x_{i+1}^2 x_{i+4}^2 \ , \ \
\right.
\\
& \qquad \left.
\sum_i x_i^3 x_{i+2} x_{i+3} \ , \ \
\sum_i x_i x_{i+2}^2 x_{i+3}^2 \ , \ \
x_0 x_1 x_2 x_3 x_4
\right> \ ,
\end{align}
where as always the subscripts on the coordinates are defined mod 5.
In particular, the well-studied Fermat quintic $\sum_i x_i^5$ and its Schoen \cite{schoen} cousin $\sum_i x_i^5 + \psi x_0 x_1 x_2 x_3 x_4$ are both illustrative examples.
In general, linear combinations of the above 6 quintics are called the Horrocks-Mumford quintic $X_{HM}$ in $\IP^4$.

We can readily compute the Molien series of $\cH(5)$ to be \cite{gap} 
\comment{
gen1 :=
[[0,0,0,0,1],[1,0,0,0,0],[0,1,0,0,0],[0,0,1,0,0],[0,0,0,1,0]];
gen2 :=
[[E(5),0,0,0,0],[0,E(5)^2,0,0,0],[0,0,E(5)^3,0,0],[0,0,0,E(5)^4,0],[0,0,0,0,1]];
heisen := Group(gen1,gen2);
irr = Irr(heisen);
###  Numbers 26, 27, 28, 29 are all 5-dim faithful representations
mol := MolienSeries(Irr(hei)[26]);
List( [ 0 .. 20 ], i -> ValueMolienSeries( mol, i ) );
}
\begin{equation}
M(z; \cH(5)) = 
\frac{ 1+z^5+21z^{10}+z^{15}+z^{20} } { (1-z^5)^5 } =
\frac{ \left[z^{10} (z^{10} + z^5 + 21) \right] \circ (z + 1/z) } { (1-z^5)^5 } 
\ ,
\end{equation}
and see that the coefficient of $z^5$ is indeed 6.
In the above, $\circ$ is functional composition.
Likewise, we can compute the Molien series, using the 5-dimensional faithful representation of $G_{HM}$, to be
\begin{equation}
M(z; G_{HM}) =
\frac{
\left[z^{35} 
(z^{35} - 2z^{30} + 2 z^{25} - z^{20} + 4 z^{15} - 4z^{10} + 12 z^{5} -9) \right] 
\circ (z + 1/z)
}{(1-z^{30})(1-z^{20})(1-z^{15})(1-z^5)^2 } 
\ .
\end{equation}
Upon developing the series, we see that the coefficient of $z^5$ is 0 so that there are no quintic invariants; the first invariant uniquely occurs at degree 10.
This indeed is why the Horrocks-Mumford quintics refer to ones composed of $\cH(5)$-invariants. We emphasize that we are considering {\it linear} invariants here. When working in projective space, here $\IP^4$, we need only consider {\it projective} invariants. For the Heisenberg group, it suffices to consider $\IZ_5 \times \IZ_5$, which is the quotient of $\cH(5)$ by its centre $\IZ_5$.
Therefore, in the literature, the HM-quintic is traditionally called the $\IZ_5 \times \IZ_5$ quintic-quotient \cite{Braun:2009mb}.

Incidentally, we notice that both Molien series have palindromic numerators.
This means that geometrically, considering them as the Hilbert series of the affine varieties $\IC^5 / \cH(5)$ and $\IC^5 / G_{HM}$ respectively, these quotients are affine (singular) Calabi-Yau \cite{stanley,Forcella:2008bb}.
This is consistent with the factor that all our explicit matrix generators, and hence all group elements, have unit determinant, and thus $\cH(5)$ and $G_{HM}$ are discrete finite subgroups of $SU(5)$. Therefore, our quotients are local Calabi-Yau 5-fold orbifolds \cite{Hanany:1998sd}.

\paragraph{Defining Equations for $X_{HM}$: }
In terms of the section of the bundle $F_{HM}$, it was shown in \cite{HM} that if $s_1 = (s_{11}, s_{12})$ and $s_2 = (s_{21}, s_{22})$ are generic sections, here written as 2-vectors because $F_{HM}$ is rank 2, then Horrocks-Mumford quintics can be written as 
\begin{equation}
s_{11} s_{22} - s_{12} s_{21} = 0
\end{equation}
and thus have nodal singularities: there are in fact 100 of them.
The form above suggests that $X_{HM}$ might be determinantal varieties. This is indeed the case \cite{lee}.
Defining the matrices
\begin{equation}
M_y(x)_{ij} := \left\{ y_{3(i-j)} x_{3(i+j)} \right\} \ ,
\qquad
L_y(z)_{ij} := \left\{ y_{i-j}z_{2i-j} \right\}
\end{equation}
for $y_j$ projective coordinates on $\IP^4_{[y_j]}$ and $z_j$ projective coordinates on $\IP^4_{[z_j]}$, we have  (note that $M_y(x)z = L_y(z)x$ with $x$ and $z$ treated as column vectors) that
\begin{equation}
\{ \det M_y(x) = 0 \} \subset \IP^4_{[x_i]} \ , \qquad
\{ \det L_y(z) = 0 \} \subset \IP^4_{[z_i]}
\end{equation}
are Horrocks-Mumford quintics.
It was shown, incidentally, that a particular blow-up of a HM quintic has (the Mellin transform of) its L-function being the unique weight 4, level 55 modular form \cite{lee}.

Moreover, in light of the sextic of Bring in \eqref{bring} and the octavic of Fricke in \eqref{frickeC} in the context of Observation \ref{360}, it is natural to consider the decimic on the Fermat Calabi-Yau, viz.,
\begin{equation}
\cD = \{
\sum_i x_i^5 = \sum_i x_i^2 = \sum_i x_i = 0
\} \subset \IP^4 \ ,
\end{equation}
of genus 16.
It would be interesting to find out how many tritangent planes (the number of bitangents, by \eqref{theta}, is $2^{15}(2^{16}-1)$) are there on $\cD$, a classical though somewhat tedious exercise which should be performed.

\subsection{Embedding into Conway's Group}
We can immediately check that of all the sporadic groups only $Co_1$, the first Conway group and $HN$, the Harada-Norton group, contain conjugacy classes whose centralizer is of order 15000.
These are class 5C of $Co_1$ and classes 5C and 5D of $HN$.
The natural question is then whether the centralizer is precisely $G_{HM}$.

Using \cite{gap}, we can readily check by direct computation \footnote{
We are grateful to Alexander Hulpke for kind help with dealing with the AtlasRep package in GAP for the simple groups of such large order.
} that this is indeed so for $Co_1$ and that the two cases for $HN$ are not.
We conclude therefore 
\begin{observation}\label{HMCo1}
The Horrocks-Mumford group $G_{HM}$ is the centralizer for exactly one conjugacy class of precisely one sporadic group: namely class 5C of $Co_1$.
\end{observation}


~\\
~\\
~\\

\section*{Acknowledgments}
We are much indebted to Professors Scott Carnahan, Chris Cummins, Igor Dolgachev, John Duncan, Gerald H\"ohn, Alexander Hulpke, Rodrigo Matias, and Simon Norton for many helpful correspondences, provision of various data, and patient help with GAP. YHH is also grateful to Dr.~Elizabeth Hunter He for painstakingly cross-checking the various tables.
YHH would like to thank the Science and Technology Facilities Council, UK, for grant ST/J00037X/1, the Chinese Ministry of Education, for a Chang-Jiang Chair Professorship at NanKai University as well as the City of Tian-Jin for a Qian-Ren Scholarship, the US NSF for grant CCF-1048082, as well as City University, London and Merton College, Oxford, for their enduring support. JM is grateful to the NSERC of Canada.

\newpage

\appendix

\section{Character Table of $G_{HM}$}\label{ap:ct}\setall
At the very end of this paper, due to its size, we present the (full ordinary linear) character table of the Horrocks-Mumford group in GAP \cite{gap} notation for reference. The notation is standard to GAP: $E(5) = \omega_5$, the primitive 5-th root of unity, ``.'' means 0, $/x$ means $1/x$ and $*x$, the complex conjugate of $x$.
We see there are 50 irreducible representations $X.1$ to $X.50$, and thus likewise 50 conjugacy classes.
The top three (unlabelled) rows of the table are respectively the class number, size and order of the conjugacy classes.

\begin{table}[h!!!]
{\scriptsize
\[\arraycolsep=1.1pt\def\arraystretch{0.5}
\begin{array}{c}
\left(
\begin{array}{rrrrrrrrrrrrrrrrrrrrrrrrr}
 1 & 1 & 1 & 1 & 1 & 1 & 1 & 1 & 1 & 1 & 1 & 1 & 1 & 1 & 1 & 1 & 1 & 1 & 1
   & 1 & 1 & 1 & 1 & 1 & 1 \\
 4 & -4 & -2 & 0 & 4 & 4 & 4 & 4 & -1 & -1 & -1 & -1 & -1 & -1 & -1 & -1 &
   -1 & -1 & 4 & -1 & -1 & -1 & -1 & 2 & -4 \\
 4 & 4 & 1 & 0 & 4 & 4 & 4 & 4 & -1 & -1 & -1 & -1 & -1 & -1 & -1 & -1 &
   -1 & -1 & 4 & -1 & -1 & -1 & -1 & 1 & 4 \\
 4 & -4 & 1 & 0 & 4 & 4 & 4 & 4 & -1 & -1 & -1 & -1 & -1 & -1 & -1 & -1 &
   -1 & -1 & 4 & -1 & -1 & -1 & -1 & -1 & -4 \\
 5 & 5 & -1 & 1 & 5 & 5 & 5 & 5 & 0 & 0 & 0 & 0 & 0 & 0 & 0 & 0 & 0 & 0 &
   5 & 0 & 0 & 0 & 0 & -1 & 5 \\
 6 & 6 & 0 & -2 & 6 & 6 & 6 & 6 & 1 & 1 & 1 & 1 & 1 & 1 & 1 & 1 & 1 & 1 &
   6 & 1 & 1 & 1 & 1 & 0 & 6 \\
 6 & -6 & 0 & 0 & 6 & 6 & 6 & 6 & 1 & 1 & 1 & 1 & 1 & 1 & 1 & 1 & 1 & 1 &
   6 & 1 & 1 & 1 & 1 & 0 & -6 \\
 20 & 4 & -4 & -4 & -5 & -5 & -5 & -5 & 5 & -5 & -5 & 5 & 0 & 5 & -5 & 0 &
   -5 & 5 & 0 & 0 & 0 & 0 & 0 & 4 & -1 \\
 24 & 0 & 0 & 0 & 24 & 24 & 24 & 24 & 4 & 4 & 4 & 4 & 4 & 4 & 4 & 4 & 4 &
   4 & -1 & -1 & -1 & -1 & -1 & 0 & 0 \\
 40 & -8 & 4 & 0 & -10 & -10 & -10 & -10 & -5 & 0 & 0 & -5 & 10 & -5 & 0 &
   10 & 0 & -5 & 0 & 0 & 0 & 0 & 0 & 4 & 2 \\
 40 & -8 & 4 & 0 & -10 & -10 & -10 & -10 & 0 & 5 & 5 & 0 & -10 & 0 & 5 &
   -10 & 5 & 0 & 0 & 0 & 0 & 0 & 0 & 4 & 2 \\
 60 & 12 & 0 & 4 & -15 & -15 & -15 & -15 & 0 & -5 & -5 & 0 & 10 & 0 & -5 &
   10 & -5 & 0 & 0 & 0 & 0 & 0 & 0 & 0 & -3 \\
 60 & 12 & 0 & 4 & -15 & -15 & -15 & -15 & 5 & 0 & 0 & 5 & -10 & 5 & 0 &
   -10 & 0 & 5 & 0 & 0 & 0 & 0 & 0 & 0 & -3 \\
 80 & -16 & -4 & 0 & -20 & -20 & -20 & -20 & -5 & 5 & 5 & -5 & 0 & -5 & 5
   & 0 & 5 & -5 & 0 & 0 & 0 & 0 & 0 & -4 & 4 \\
 80 & 16 & -4 & 0 & -20 & -20 & -20 & -20 & -5 & 5 & 5 & -5 & 0 & -5 & 5 &
   0 & 5 & -5 & 0 & 0 & 0 & 0 & 0 & 4 & -4 \\
 96 & 0 & 0 & 0 & 96 & 96 & 96 & 96 & -4 & -4 & -4 & -4 & -4 & -4 & -4 &
   -4 & -4 & -4 & -4 & 1 & 1 & 1 & 1 & 0 & 0 \\
 100 & 20 & 4 & -4 & -25 & -25 & -25 & -25 & 0 & 0 & 0 & 0 & 0 & 0 & 0 & 0
   & 0 & 0 & 0 & 0 & 0 & 0 & 0 & -4 & -5 \\
 120 & -24 & 0 & 0 & -30 & -30 & -30 & -30 & 5 & -5 & -5 & 5 & 0 & 5 & -5
   & 0 & -5 & 5 & 0 & 0 & 0 & 0 & 0 & 0 & 6 \\
\end{array}
\right. \ldots
\\
\\
\\
\ldots
\left.
\begin{array}{rrrrrrrrrrrrrrrrrrrrrrrrr}
1 & 1 & 1 & 1 & 1 & 1 & 1 & 1 & 1 & 1 & 1 & 1 & 1 & 1 & 1 & 1 & 1 & 1 & 1
   & 1 & 1 & 1 & 1 & 1 & 1 \\
 -4 & -4 & -4 & 1 & 1 & 1 & 1 & 1 & 1 & 1 & 1 & 1 & 1 & -2 & -2 & -2 & -2
   & 0 & 0 & 0 & 0 & 2 & 2 & 2 & 2 \\
 4 & 4 & 4 & -1 & -1 & -1 & -1 & -1 & -1 & -1 & -1 & -1 & -1 & 1 & 1 & 1 &
   1 & 0 & 0 & 0 & 0 & 1 & 1 & 1 & 1 \\
 -4 & -4 & -4 & 1 & 1 & 1 & 1 & 1 & 1 & 1 & 1 & 1 & 1 & 1 & 1 & 1 & 1 & 0
   & 0 & 0 & 0 & -1 & -1 & -1 & -1 \\
 5 & 5 & 5 & 0 & 0 & 0 & 0 & 0 & 0 & 0 & 0 & 0 & 0 & -1 & -1 & -1 & -1 & 1
   & 1 & 1 & 1 & -1 & -1 & -1 & -1 \\
 6 & 6 & 6 & 1 & 1 & 1 & 1 & 1 & 1 & 1 & 1 & 1 & 1 & 0 & 0 & 0 & 0 & -2 &
   -2 & -2 & -2 & 0 & 0 & 0 & 0 \\
 -6 & -6 & -6 & -1 & -1 & -1 & -1 & -1 & -1 & -1 & -1 & -1 & -1 & 0 & 0 &
   0 & 0 & 0 & 0 & 0 & 0 & 0 & 0 & 0 & 0 \\
 -1 & -1 & -1 & -1 & 4 & -1 & -1 & -1 & -1 & -1 & -1 & -1 & 4 & 1 & 1 & 1
   & 1 & 1 & 1 & 1 & 1 & -1 & -1 & -1 & -1 \\
 0 & 0 & 0 & 0 & 0 & 0 & 0 & 0 & 0 & 0 & 0 & 0 & 0 & 0 & 0 & 0 & 0 & 0 & 0
   & 0 & 0 & 0 & 0 & 0 & 0 \\
 2 & 2 & 2 & 2 & 2 & -3 & 2 & 2 & -3 & -3 & -3 & 2 & 2 & -1 & -1 & -1 & -1
   & 0 & 0 & 0 & 0 & -1 & -1 & -1 & -1 \\
 2 & 2 & 2 & -3 & 2 & 2 & -3 & -3 & 2 & 2 & 2 & -3 & 2 & -1 & -1 & -1 & -1
   & 0 & 0 & 0 & 0 & -1 & -1 & -1 & -1 \\
 -3 & -3 & -3 & -3 & 2 & 2 & -3 & -3 & 2 & 2 & 2 & -3 & 2 & 0 & 0 & 0 & 0
   & -1 & -1 & -1 & -1 & 0 & 0 & 0 & 0 \\
 -3 & -3 & -3 & 2 & 2 & -3 & 2 & 2 & -3 & -3 & -3 & 2 & 2 & 0 & 0 & 0 & 0
   & -1 & -1 & -1 & -1 & 0 & 0 & 0 & 0 \\
 4 & 4 & 4 & -1 & 4 & -1 & -1 & -1 & -1 & -1 & -1 & -1 & 4 & 1 & 1 & 1 & 1
   & 0 & 0 & 0 & 0 & 1 & 1 & 1 & 1 \\
 -4 & -4 & -4 & 1 & -4 & 1 & 1 & 1 & 1 & 1 & 1 & 1 & -4 & 1 & 1 & 1 & 1 &
   0 & 0 & 0 & 0 & -1 & -1 & -1 & -1 \\
 0 & 0 & 0 & 0 & 0 & 0 & 0 & 0 & 0 & 0 & 0 & 0 & 0 & 0 & 0 & 0 & 0 & 0 & 0
   & 0 & 0 & 0 & 0 & 0 & 0 \\
 -5 & -5 & -5 & 0 & 0 & 0 & 0 & 0 & 0 & 0 & 0 & 0 & 0 & -1 & -1 & -1 & -1
   & 1 & 1 & 1 & 1 & 1 & 1 & 1 & 1 \\
 6 & 6 & 6 & 1 & -4 & 1 & 1 & 1 & 1 & 1 & 1 & 1 & -4 & 0 & 0 & 0 & 0 & 0 &
   0 & 0 & 0 & 0 & 0 & 0 & 0 \\
\end{array}
\right)
\end{array}\]
}
\caption{\sf
The rational character table for $G_{HM}$.
\label{t:rct}
}
\end{table}

We can combine the irreducible representations (rows) of the above ordinary character table in the following groups
    \{ 1 \}, \
    \{ 2, 3 \}, \
    \{ 6 \}, \
    \{ 7 \}, \
    \{ 8 \}, \
    \{ 4, 5 \}, \
    \{ 13 \}, \
    \{ 9, 10, 11, 12 \}, \
    \{ 38 \}, \
    \{ 14, 16, 18, 19 \}, \
    \{ 15, 17, 20, 21 \}, \
    \{ 22, 24, 26, 28 \}, \
    \{ 23, 25, 27, 29 \}, \
    \{ 30, 34, 35, 37 \}, \
    \{ 31, 32, 33, 36 \}, \
    \{ 39, 40, 41, 42 \}, \
    \{ 43, 44, 45, 46 \}, \
    \{ 47, 48, 49, 50 \}
to produce a rational character table wherein the Galois conjugates conspire to cancel.
This is done using MAGMA \cite{magma} and presented in Table \ref{t:rct}.

For completeness, we also present the Brauer Character Table in $\IF_5$, the choice of this latter field being obvious from the definition of $G_{HM}$:
\[
\arraycolsep=1pt\def\arraystretch{0.6}
\left(
\begin{array}{rrrr}
 1 & 1 & 1 & 1 \\
 2 & -2 & -1 & 0 \\
 3 & 3 & 0 & -1 \\
 4 & -4 & 1 & 0 \\
 5 & 5 & -1 & 1 \\
\end{array}
\right.
\quad
\ldots
\quad
\mbox{{\Huge $0_{{\mbox{\normalsize $5 \times 19$}}}$}}
\quad
\ldots
\quad
\left.
\begin{array}{r}
 1 \\
 1 \\
 0 \\
 -1 \\
 -1 \\
\end{array}
\quad
\mbox{{\Huge $0_{{\mbox{\normalsize $5 \times 26$}}}$}}
\right)_{5 \times 50}
\]
We see that there are only 5 irreducible representations.

\subsection{Details of Construction}\label{s:HMdetail}
A little more detail about the Horrocks-Mumford construction is nicely summarized in \cite{lee}, which we recapitulate briefly here.
Take $V = \IC^5$ (with standard basis $e_{i=1,\ldots,5}$) so that $\IP^4 = \IP(V)$. We have the Koszul sequence
\begin{equation}
0 \to \cO \to V\otimes\cO(1) \to \wedge^2 V \otimes\cO(2)
\wedge^3 V \otimes\cO(3) \to \wedge^4 V \otimes\cO(4) \to \cO(5) \to 0
\end{equation}
where all sheafs are on $\IP^4$.
The kernel of the first map $\cO \to V\otimes\cO(1)$ is none other than the tangent bundle $T$ of $\IP^4$ and more generally $im \left( \cO(i) \otimes \wedge^i V \right) \simeq \wedge^i T$.
We hence have the sequence
\begin{equation}
\wedge^2 V \otimes\cO(2) 
\stackrel{p_0}{\longrightarrow}
\wedge^2 T
\stackrel{q_0}{\longrightarrow} 
\wedge^3 V \otimes\cO(3) 
\end{equation}
with $p_0$ surjective and $q_0$ injective.
Horrocks-Mumford then defines the maps
\begin{align}
\nn
p & : V \otimes \cO(2) 
  \stackrel{{\scriptsize \left(\begin{array}{c}f^+\\f^-\end{array}\right)}\otimes\cO(2)}{\xrightarrow{\hspace*{3cm}}} 
(\wedge^2 V)^{\oplus 2} \otimes \cO(2)
  \stackrel{p_0^{\oplus 2}}{\longrightarrow} 
(\wedge^2 T)^{\oplus 2}
\ ;\\
\nn
q & : (\wedge^2V)^{\oplus 2} 
  \stackrel{q_0^{\oplus 2}}{\longrightarrow} 
(\wedge^3V)^{\oplus 2} \otimes \cO(3)
  \stackrel{{\scriptsize \left( -f^{-*}, \ f^{+*} \right)}\otimes\cO(3)}{\xrightarrow{\hspace*{3cm}}} 
V^* \otimes \cO(3)
\ ;
\\
\nn & \qquad \mbox{with}
\\
& 
\begin{array}{cccc}
  & V & \longrightarrow & \wedge^2 V \\
f^+ : &\sum\limits_{i=1}^5 v_i e_i &\longrightarrow&
      \sum\limits_{i=1}^5 v_i e_{i+2} \wedge e_{i+3} \ ,
\\
f^- : &\sum\limits_{i=1}^5 v_i e_i &\longrightarrow&
      \sum\limits_{i=1}^5 v_i e_{i+1} \wedge e_{i+4} \ .
\end{array}
\end{align}
Clearly $q \circ p = 0$ and we have the complex
\begin{equation}
V \otimes \cO(2) 
\stackrel{p}{\longrightarrow}
(\wedge^2 T)^{\oplus 2}
\stackrel{q}{\longrightarrow}
V^* \otimes \cO(3) \ . 
\end{equation}
Recalling that $V = \IC^5$, we retrieve \eqref{HMseq} and \eqref{HMdef}.


\topmargin       -1.0in  
\textheight      10in  
\thispagestyle{empty}
\[
\begin{array}{ll}
\begin{array}{l}
\includegraphics[trim=0mm 0mm 0mm 0mm, width=9.3in,angle=270]{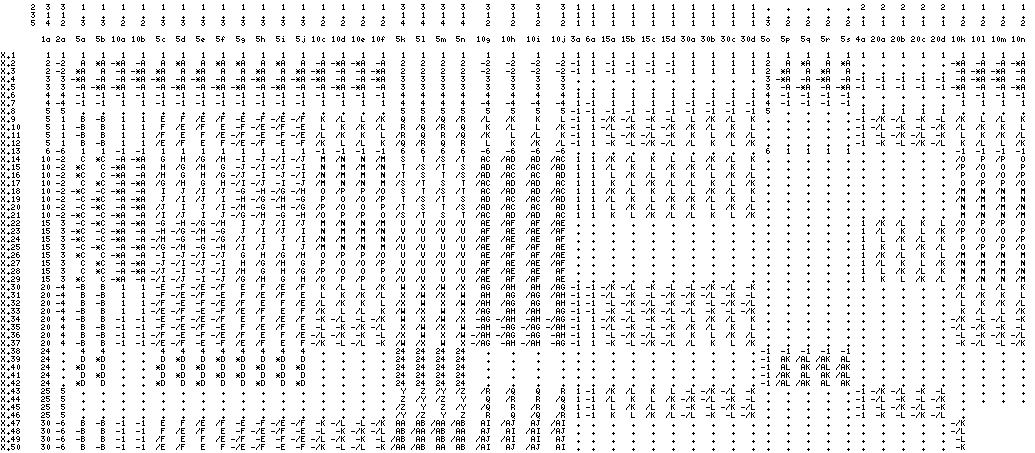}
\end{array}
&
\begin{array}{l}
\includegraphics[trim=0mm 0mm 0mm 0mm, width=2.4in]{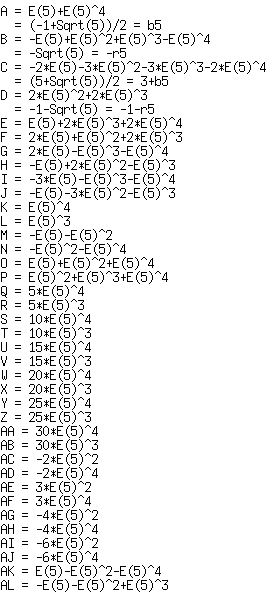}
\end{array}
\end{array}
\]

\end{document}